\documentclass[10pt]{article}

\usepackage{amssymb}
\usepackage{amsmath}
\usepackage{amsthm}
\usepackage{graphicx}
\usepackage[per-mode=symbol]{siunitx}
\usepackage{mathtools}
\usepackage{tikz}
\usetikzlibrary{arrows,calc}
\usetikzlibrary{shapes,arrows,shadows,arrows.meta,decorations.pathreplacing,angles,quotes}
\usepackage{bm}
\usepackage{booktabs}
\usepackage{pgfplots}
\pgfplotsset{compat=1.15}
\usepackage{float}
\usepackage[section]{placeins}

\usepackage{lineno,hyperref}
\modulolinenumbers[5]
\usepackage{subcaption} 
\usepackage{url} 

\newcommand{\txt}[1]{\text{#1}} 

\newcommand{\mat}[1]{{\mathbf{#1}}}
\newcommand{\cur}[1]{{\mathcal{#1}}}

\renewcommand{\d}{\,\text{d}}

\newcommand{\TOL}{\,\texttt{TOL} }

\newcommand{\Jpde}{\ensuremath{\cur J_{\textrm{P\!\,D\!\,E}}}}
\newcommand{\J}{\ensuremath{\cur J}}

\newcommand{\atm}{{\text{\,atm}}}

\newcommand{\Had}{\ensuremath{H_{\text{ad}}}}
\newcommand{\etaad}{\ensuremath{\eta_{\text{ad}}}}
\newcommand{\ds}{\displaystyle}
\newcommand{\R}{{\mathbb R}}

\pgfplotsset{every axis legend/.append style={
at={(1.02,1)},
anchor=north west}}

\tikzset{
   compressor/.pic={
      \draw[line width = 1pt, fill=white] (0,0) circle (0.2);
      \draw[line width = 1pt] (10:0.2) -- (130:0.2);
      \draw[line width = 1pt] (-10:0.2) -- (-130:0.2);
   },
   valve/.pic={
      \draw[line width = 1pt, fill=white] (-0.2,-0.15) -- (0.2,0.15) -- (0.2,-0.15) -- (-0.2,0.15) -- cycle;
   },
   controlValve/.pic={
      \draw[line width = 1pt, fill=white] (-0.2,-0.15) rectangle (0.2,0.15) ;
      \draw[line width = 1pt] (-0.2,-0.15) -- (0.2,0) -- (-0.2,0.15);
   },
   source/.pic={
      \draw[line width = 1pt, fill=white] (0,0) circle (0.15);
   },
   sink/.pic={
      \draw[line width = 1pt, fill=white] (-0.2,-0.2) -- (0.2,-0.2) -- (0.2,0.1) -- (0,0.3) -- (-0.2,0.1) -- cycle;
      \draw[line width = 1pt] (-0.2,0.1) -- (0.2,0.1);
   }
}

\def\colALG{black!20}
\def\colLIN{black!40}
\def\colNL{black!60}

\graphicspath{{figures/}}

\author{P.~Domschke, O.~Kolb, J.~Lang}
\title{Fast and Reliable Transient Simulation and Continuous Optimization of Large-Scale Gas Networks}
\author{
Pia Domschke\\
{\small \it Frankfurt School of Finance \& Management, Finance Department,} \\
{\small Adickesallee 32-34, 60322 Frankfurt am Main, Germany}\\
{\small p.domschke@fs.de} \\ \\
Oliver Kolb\\
{\small \it University of Mannheim, Department of Mathematics,} \\
{\small A5, 6, 68131 Mannheim, Germany}\\
{\small kolb@uni-mannheim.de} \\ \\
Jens Lang\footnote{corresponding author, ORCID: 0000-0003-4603-6554}\\
{\small \it Technical University of Darmstadt, Department of Mathematics,} \\
{\small Dolivostra{\ss}e 15, 64293 Darmstadt, Germany}\\
{\small lang@mathematik.tu-darmstadt.de}}
\begin{document}
\maketitle

\begin{abstract}
We are concerned with the simulation and optimization of large-scale gas pipeline
systems in an error-controlled environment. The gas flow dynamics is locally
approximated by sufficiently accurate physical models taken from a hierarchy of
decreasing complexity and varying over time. Feasible work regions of compressor
stations consisting of several turbo compressors are included by semiconvex
approximations of aggregated characteristic fields. A discrete adjoint approach
within a first-discretize-then-optimize strategy is proposed and a sequential
quadratic programming with an active set strategy is applied to solve the
nonlinear constrained optimization problems resulting from a validation of nominations.
The method proposed here accelerates the computation of near-term forecasts of
sudden changes in the gas management and allows for an economic control of intra-day
gas flow schedules in large networks. Case studies for real gas pipeline systems
show the remarkable performance of the new method.
\end{abstract}



\newpage

\section{Introduction}
The ongoing replacement of traditional energy production by coal fired and nuclear plants
with gas consuming facilities has rapidly increased the role of natural gas transport through
large networks. The security of energy supply and the development of clean energy to meet
environmental demands have generated a significant increase in gas consumption for electric
power stations in the last decade. In these days, natural gas is considered as a bridging
combustible resource on the way towards a future energy mix mainly based on low-carbon and
regenerative energy \cite{EA2019}. The seasonally fluctuating disposability of wind and solar resources
causes a growing variability in electricity production and hence also in the demands of gas
transportation by pipelines. The resulting intra-day oscillations in demand for natural gas
leads to new challenges for computer based modelling and control of gas pipeline operations
with an increasing focus on short-term transient dynamics of gas flow. Operators have to
responsively control loads to realize a reliable operational management for both gas and
electricity delivery systems. The challenging new conditions demand advanced decision tools
based on transient nonlinear optimization taking into account serious operating restrictions.

In this paper, we propose a novel computational approach for the reliable transient simulation
and continuous optimization for gas pipeline flow through networks. The operation of compressor
stations and flow schedules are determined in such a way that operating limits of compressors
and pressure constraints inside the pipes are satisfied. Exemplarily, we investigate the important
task of safely driving a stationary running system from an initial condition to a target state
defined by shifted
gas nominations at the entrances or delivery points of the network. To be usable in real-time
gas management, we have designed our methods to meet user-defined accuracies while keeping the
computing time for large-scale gas networks at a moderate level. The proposed method is fast
enough to allow the application of uncertainty quantification by Monte Carlo or other stochastic
simulations for risk analysis and reliability assessment of gas delivery.

The highly nonlinear partial differential equations representing the constraints of gas
flow dynamics are locally approximated by a certain
flow model taken from a whole hierarchy of decreasing complexity. Error estimates are computed
from sensitivity equations and used as the basis for a straightforward criterion to verify that
the model chosen provides a physically meaningful representation of the local gas flow in a pipe.
Similar estimators are exploited to set up a fully adaptive spatial resolution in each
individual pipe and variable time steps in a finite volume setting, which greatly reduces the
size of the resulting discretized system. These strategies have been already proven to work
reliably for networks of academic size \cite{DomschkeDuaStolwijkEtAl2018,DomschkeKolbLang2015}.

Multiple compressor stations consisting of several compressor machines often connected through
intricate topological designs compensate the pressure drop due to the inherent friction or
user's excessive withdrawal of natural gas. Additionally, complex constraint envelopes delimit
the allowed operating states of the compressors, making a generalized modelling and managing
of such stations a challenging task of its own. To enable an efficient usage within a nonlinear
optimization, we extend the approximate convex decomposition developed in
\cite{HillerWalther2017,LienAmato2006} to a semiconvex setting. Critically, a semiconvex
approximation of the characteristic work diagram for a compressor station enables a
new representation of compressor constraints and can be used for a fast validation
of nominations.

The development of an optimization-based, automated decision system requires the calculation
of gradient information to be completed repeatedly and reliably in a short time in order to
provide a decision support that is usable in a real-time gas management for large-scale
pipeline systems. We propose to use a discrete adjoint approach within a first-discretize-then-optimize
strategy applied to a possible splitting of the overall time horizon into smaller subintervals.
A sequential quadratic programming with an active set strategy is then applied to solve the
resulting nonlinear constrained optimization problems \cite{Spellucci1998,Spellucci1998a}.

There are many academic and industrial studies available in the literature for the
gas transport through a single pipe or even whole gas networks. In the last fifteen years, an increasing
effort for the optimization of transient processes has led to a couple of publications that consider
gas pipeline systems of moderate size \cite{BurlacuEgger2019,DomschkeGeislerKolbEtAl2011,DomschkeKolbLang2015,EhrhardtSteinbach2005,
HanteLeugeringMartinEtAl2017,Herty2007,MakHentenryck2019,Pfetsch2014}. Very recently, a linearization procedure that enables the use of closed-form solutions for the gas transport equations
was proposed in \cite{BeylinRudkevichZlotnik2020} and applied to middle-sized networks. Significant progress has been made in developing sophisticated models for compressor stations, which is a key
for more realistic treatment of an adjusted pressure increase under certain operating limits
\cite{BeylinRudkevichZlotnik2020,RoseSchmidtSteinbachEtAl2016,WaltherHiller2017}. The gas flow in larger pipeline systems may
be regulated by control valves and groups of compressors could be inactive in certain time
intervals. Optimization of such systems leads to transient mixed-integer nonlinear problems. Although
significant progress has been made to tackle these challenging problems in recent years, they still
remain unsolved for real-life applications with sudden changes in gas supply and demand, which
requires detailed transport models and reliable numerical solutions \cite{BurlacuEgger2019,DomschkeGeislerKolbEtAl2011,HahnLeyfferZavala2017,MahlkeMartinMoritz2010}.
Model order reduction is an appropriate technique to reduce the computational cost for
simulating large-scale gas networks. Surrogates
based on radial functions and proper orthogonal decomposition for quasi-static approximations of the
gas flow are compared in \cite{GrundelHornungKlaassenBennerClees2013}. However, the construction
of such models in the case of highly transient flow with nonlinear compressors is still very
ambitious. Recently, the tracking of internal space-time flow and pressure profiles using
pressure measurements at nodal junctions and optimization-based state estimates has been investigated in \cite{JalvingZavala2018} for simplified compressor models and networks of moderate size.
Graph-based modeling abstractions with application to optimal control of connected gas pipelines in series
is proposed in \cite{JalvingCaoZavala2019}. The challenging
task of maximizing the economic welfare of gas users for slow transients that allow a simplification
of Euler's equations is the objective in \cite{ZlotnikSundarRudkevichBeylinLi2019}. To the best
of our knowledge, in all these approaches there is no rigorous control of the discretization errors
caused by the space-time resolution and model selection.

We apply our novel computation and optimization approach to real networks taken from the GasLib library
\cite{SchmidtAsmannBurlacuEtAl2017}. Fully model-space-time-adaptive simulations provide
accuracy bounds for target functionals and can be viewed as an efficient way to automatically and
safely reduce the order of the gas network model. They can be computed in the range of seconds for several
hundreds of edges including pipes, valves, and compressor stations. Making use of this
acceleration and calling the simulation tool within a nonlinear continuous optimization solver to
determine gradient information reduces the computing time to validate nominations by several orders
of magnitude compared to previously developed methods. The approach may be promising to include
binary decision variables and to perform reliability studies using methods from uncertainty
quantification in future extensions.

The paper is structured as follows. In Section 2, we outline our gas network modelling including
boundary and coupling conditions, modelling of gas flow through pipes, compressor stations, and
valves. Our numerical schemes for simulation and optimization are described in Section 3. In
Section 4 and 5, we present results for certain case studies including nomination validation.
Finally, a short summary and a outlook
are given in Section 6.

\newpage

\section{Gas Network Modelling} \label{sec:modelling}
We model gas supply networks as a directed graph $\cur G = (\cur J, \cur V)$ with arcs~$\cur J$ and vertices~$\cur V$ (nodes, branching points).
The set of arcs~$\cur J$ contains pipes, compressor stations, simple valves and control valves. The gas dynamics within the pipes is described by a hierarchy of models of different complexity ranging from simple algebraic equations to a hyperbolic system of PDEs (see below and compare to~\cite{DomschkeKolbLang2015}). For each arc in the network that is modelled by a PDE, we consider an interval $[x_j^a, x_j^b]$ with $x_j^a < x_j^b$, $j \in \Jpde$, where $\Jpde\subseteq \J$ contains all arcs modelled with PDEs. All other network components, including coupling and boundary conditions at the nodes, are described by algebraic equations, where we also use spatial coordinates $x_j^a$ and $x_j^b$ to describe states at the beginning (tail) and end (head) of an arc $j \in \cur J \setminus \Jpde$ or denote states at the beginning/end of the arc with subscripts ``in''/``out''.
Note that one has to specify adequate initial, coupling and boundary conditions for a complete problem description of the gas network.

The flow through the network is described by the state variables $q(x,t)$ and $p(x,t)$. Here,
$q$ denotes the flow rate (or flux) in \si{\meter\cubed\per\second} under standard conditions, i.e., pressure of \SI{1}{\atm} and temperature of $\SI{0}{\celsius}=\SI{273.15}{\kelvin}$, and $p$ represents the pressure in \si{\pascal}.\\

\noindent\textit{Coupling and Boundary Conditions}.
For any node $v\in \cur{V}$, we denote by $\delta_v^-$ the ingoing arcs of $v$ and the outgoing arcs by $\delta_v^+$. Respecting conservation of mass, the entire flux going into any node has to be equal to the sum of fluxes going out of that node:
\begin{equation}
 \sum\limits_{j\in\delta_v^+} q(x_j^a,t) - \sum\limits_{j\in\delta_v^-} q(x_j^b,t) = q(v,t) \quad \forall t>0 \label{equ:conservationOfMass_generalized}
\end{equation}
where $q(v,t)$ is an auxiliary variable to model the flow rate of a feed-in or demand at node $v$ (see below). Otherwise, it is set to zero. In addition to conservation of mass, it is widely used to claim the equality of pressure at all nodes $v \in \cur V$, that is,
\begin{equation}
p(x_j^a,t) = p(v,t) \quad \forall j\in
\delta_v^+, \quad p(x_i^b,t) = p(v,t) \quad \forall i\in \delta_v^-,
\label{equ:equalityOfPressure}
\end{equation}
with an auxiliary variable $p(v,t)$ for the pressure at the node $v$.
For coupling and boundary nodes, one has to pose one more condition,
which we consider of the form
\begin{equation}
e(p(v,t),q(v,t)) = 0\,.
\label{equ:node-condition}
\end{equation}
The latter allows for prescribing pressure or flux profiles as boundary condition or feed-in/demand.
For example, the pressure is often described at the inflow boundary and certain flux profiles are
requested by the gas consumers at the outflow boundary.
Other coupling conditions like conservation of energy or entropy for isothermal, isentropic or polytropic gas flow, respectively, are discussed in \cite{Egger2018,LangMindt2018,MindtLangDomschke2019}.\\

\noindent\textit{Modelling of Gas Flow through Pipes}.
Typically, gas pipeline systems are buried underground and hence temperature differences
between a pipe segment and the ground can be neglected in practice. It is therefore standard
to consider an isothermal process without a conservation law for the energy.
As most complex model for the gas dynamics in the pipes, we take

\begin{itemize}
\item (M3), the scaled isothermal Euler equations
\begin{align}
 u_t + (\mat{A}u + B(u))_x = g(u)\,, \label{equ:nonlinearModel}
\end{align}
with
\begin{align*}
 u &= \begin{pmatrix}p\\q\end{pmatrix}\!,
 &\mat{A} &= \begin{pmatrix}0 & \ds\frac{\rho_0 c^2}{A}\\[2mm] \ds \frac{A}{\rho_0} & 0\end{pmatrix}\!,
 &B(u) &=  \begin{pmatrix}0\\[2mm] \ds\frac{\rho_0 c^2q^2}{Ap} \end{pmatrix}\!,
 &g(u) &=  \begin{pmatrix}0\\[2mm] \ds -\frac{\lambda \rho_0 c^2 |q|q}{2dAp}\end{pmatrix}\!,
\end{align*}
\end{itemize}
together with the equation of state for real gases, $p=\rho zRT$, with the density $\rho$, the compressibility factor $z\in(0,1)$ \cite{KochHillerPfetschEtAl2015}, the temperature $T$, and the specific gas constant $R$.
Here, $c=\sqrt{p/\rho}$~denotes the speed of sound, $\lambda$~the friction coefficient, $d$~the diameter, $A$~the cross-sectional area of the pipe, and $\rho_0$ the density under standard conditions.
Neglecting the nonlinear term $B(u)$ on the left hand side leads to
\begin{itemize}
\item (M2), the semilinear isothermal Euler equations
\begin{align}
 u_t + \mat{A}u_x = g(u)\,.
\end{align}
\end{itemize}
As third and most simplest model, we consider
\begin{itemize}
\item (M1), the (quasi-)stationary semilinear isothermal Euler equations
\begin{align}
 \mat{A}u_x = g(u)\,,
\end{align}
\end{itemize}
which can be solved analytically,
\begin{subequations}
\begin{eqnarray}
    q    &=& \text{const.} \,,\\
    p_\txt{out} &=& \sqrt{p_\txt{in}^2 - \frac{\lambda \rho_0^2 c^2 L}{dA^2}|q|q}.\label{equ:algebraicModel-b}
\end{eqnarray} \label{equ:algebraicModel}
\end{subequations}%
Here, $L$ denotes the length of the pipe. The flow rate $q$ is now only represented by the
pressure difference in a pipe - the so-called Weymouth equation. It is the standard model to
capture long-term planning behaviour of gas networks in the quasi-stationary regime.

To summarize, the three models
(M3)-(M2)-(M1) forms a hierarchy of gas transport models with decreasing fidelity. For later use, it
is important to mention that they are characterized by an additive structure, i.e., they are related by
adding or subtracting certain differential terms. Given a reliable accuracy control, they can be used in
each individual pipe adopted to the local flow behaviour without any difficulty, since the state vector
$u=(p,q)$ remains unchanged and hence two pipes can be easily interconnected by the continuity of $p$ and $q$. \\

\noindent\textit{Modelling of Compressor Stations}.
Compressor stations consist of at least one turbo compressor, where typically various ways of controlling the entire station are available (parallel, series, subsets). Thus a direct approach of modelling a compressor station leads to mixed integer control problems
\cite{RoseSchmidtSteinbachEtAl2016}. As a simplification, an outer linear approximation of the feasible states of a compressor station is proposed in \cite{WaltherHiller2017}. This results in a characteristic diagram given by a polyhedral model in the ($Q$,$\Had$)-space, where $Q=q\rho_0/\rho$ is the volumetric flow rate in \si{\meter\cubed\per\second} and $\Had$ is the adiabatic head of the compressor defined by the
gas compression from $p_\txt{in}$ to $p_\txt{out}$:
\begin{equation}
\Had = zTR_s\frac{\kappa}{\kappa-1}\left(\left(\frac{p_\txt{out}}{p_\txt{in}}\right)^{\frac{\kappa-1}{\kappa}}
-1\right)\,.
\end{equation}
Here, $\kappa$ is the isentropic exponent and $R_s:=R/m$ is the specific gas constant with $m$ being the averaged molar mass of the gas mixture.

We make use of such compressor station models to be able to apply methods from continuous nonlinear optimization for optimal control tasks. Therefore, $\Had$ is used as time-dependent control variable, while the underlying polyhedral model is incorporated as state constraint within the optimization procedures.
The power $P$ that is needed for the compression process is given by $P=\rho_\txt{in} Q \Had/\etaad$,
where $\etaad\in[0,1]$ denotes the adiabatic efficiency depending on $Q$, the compressor speed,
and other compressor-specific parameters \cite{RoseSchmidtSteinbachEtAl2016,WaltherHiller2017}.
For later use, we denote by $\J_{cp}$ the set of all compressor stations.\\

\noindent\textit{Modelling of Valves}.
Valves are used to regulate the flow in gas networks. The equations describing an open simple valve are
\begin{equation}
q_\txt{in} = q_\txt{out}\,,\quad p_\txt{in} = p_\txt{out}\,.
\end{equation}
A closed simple valve is described by $q_\txt{in} = q_\txt{out} = 0\,.$
A special class of valves are control valves
that reduce the gas pressure by a certain amount. The corresponding equations
are given by
\begin{equation}
 p_\txt{in} - p_\txt{out} = \triangle p\,,\quad q_\txt{in} = q_\txt{out}\,,
\end{equation}
with a possibly time-dependent control variable $\triangle p=\triangle p(t)$.

\section{Numerical Schemes for Simulation and Optimization} \label{sec:numericalSchemes}
With the different models described above and all controls $(\Had(t),\triangle p(t))$ given, we can solve the whole network as a system of differential-algebraic equations using adequate initial,
coupling, and boundary conditions. Since the gas transport through a complex network may be
very dynamic and thus changes both in space and time, an automatic control of the accuracy of the simulation is mandatory. In addition, a model adaptation, i.e., switching between models (M3), (M2), and (M1) in an appropriate way, has proven to be very useful in order
to further reduce computational costs. The main idea is to use the most complex model (M3) only when necessary and to refine spatial and temporal discretizations only where needed. A complete
description of the overall strategy to efficiently control model and discretization errors up to a user-given tolerance is given in \cite{DomschkeDuaStolwijkEtAl2018,DomschkeKolbLang2015}. In what follows, we will give a brief overview on the main ingredients.\\

\subsection{Adaptive Network Simulation}
For the discretization of the hyperbolic PDEs in the pipes, we apply an implicit box scheme~\cite{KolbLangBales2010}, which reads for model (M3)
\begin{align}
\frac{u^l_{i-1}+u_i^l}{2} =&\;\frac{u^{l-1}_{i-1}+u_i^{l-1}}{2} -
\frac{\triangle t}{\triangle x} \left( \mat{A}u^l_{i}+B(u^l_{i}) - \mat{A}u^{l}_{i-1}-B(u^{l}_{i-1}) \right)
 + \frac{\triangle t}{2} \left( g(u^l_{i-1}) + g(u^l_{i}) \right).
\end{align}
It uses step sizes $\triangle x$ and $\triangle t$ in space and time, respectively, and forms a system of
nonlinear equations for the approximations $u_i^l\approx u(x_i,t_l)$. The scheme is closely related to
the finite difference method proposed in \cite{Kiuchi1994}, where the midpoint rule is used to approximate
the source term $g(u)$ -- a fact we were not aware of until now. Explicit methods as recently presented in
\cite{DyachenkoZlotnikKorotkevichChertkov2017,GyryaZlotnik2019} for large-scale natural gas pipeline networks
can be an efficient alternative if the dynamics is not highly spread over the whole network. However, optimal
control of compressor stations typically changes the dynamics close to the compressors in an unpredictable
manner, which would force explicit integrators to apply very small time steps due to the well-known CFL condition.
This effect can be even strengthened if the source term $g(u)$ becomes dominant in the turbulent region. So,
we have decided to use an implicit scheme in order to choose an appropriate time step with respect to
accuracy requirements only.

Dropping terms with $B$ on the
right-hand side yields a discretization for model (M2). No discretization for (M3) is necessary, since
it can be solved analytically. The scheme is convergent of order two in space and
order one in time. It is conservative and stable under mild conditions \cite{KolbLangBales2010}. Each pipe
of the network is discretized by using one of the three gas transport models described above and an individual
step size in space, where both can vary over time.
The global time steps taken by this scheme are also used for all algebraic model equations including model (M1).
The adaptation of the discretization parameters $h:=(\triangle x,\triangle t,(Mk)_{k\in\{1,2,3\}}$) is realized in
a successive process of the classical loop
\begin{equation}\label{semr_loop}
SOLVE \rightarrow ESTIMATE \rightarrow MARK \rightarrow REFINE\,.
\end{equation}
Given all discretization parameters collected in $h$, the model equations are solved and error
estimators are derived. These
estimators are then used to mark pipes for spatial refinement or model enhancement and to adopt the global
time step. The error estimators considered are based upon the discretized model equations and are supposed to measure the influence of the currently applied models and discretizations in each
individual pipe on a generic user-defined output functional
\begin{align}
 M(u) = \int_Q \!\! N(u) \d(x,t)
 + \sum_{v \in\cur V}\int_0^T \!\! N_{v}(u) \d t + \sum_{i\in \cur J\setminus \Jpde} \!\! \int_0^{T} \!\! N_i(u_i)\d t\,. \label{equ:functional-M-general-form}
\end{align}
Here, $Q = \Omega \times (0,T)$ with $\Omega = \bigcup_{j\in\cur \Jpde} [x_j^a,x_j^b]$ and the vector $u_i$ is defined as $u_i(t) = (p(x_i^a,t), q(x_i^a,t), p(x_i^b,t), q(x_i^b,t))^T$ for all arcs $i \in \cur{J}\setminus \Jpde$ that are modelled by algebraic equations. The functions
$N(u)$, $N_v(u)$, and $N_i(u_i)$ define tracking-type costs on the respective sets
($\Omega$, nodes, algebraic arcs) in the whole time interval $(0,T)$.
In the spirit of dual weighted residual methods, the applied error estimators are computed via sensitivity information coming from adjoint equations of the discretized model equations~\cite{DomschkeKolbLang2015}. Here, we will demonstrate their computation for the case that model (M2) is applied in all pipes. The adjoint equations of model (M2) with respect to $M(p,q)$ are given by
\begin{align}
\label{adjoint_equations}
\begin{split}
 \xi_{1_t} + \frac{A}{\rho_0} \xi_{2_x} &
   = -\frac{\lambda \rho_0 c^2}{2 D A} \frac{|q| q}{p^2} \xi_2 - N_p(p,q), \\
 \xi_{2_t} + \frac{\rho_0 c^2}{A} \xi_{1_x} &
   = \frac{\lambda \rho_0 c^2}{D A} \frac{|q|}{p} \xi_2 - N_q(p,q),
\end{split}
\end{align}
together with appropriate end, coupling, and node conditions, where
also the functions $N_v$ and $N_i$ appear \cite{Domschke2011,DomschkeKolbLang2015}.
The solution $\xi = (\xi_1,\xi_2)^T$ of the adjoint equations consists of the adjoint pressure
and flow rate of the semilinear model (M2) with respect to the functional~$M(u)$.
Let now $u=(p,q)^T$ be the solution of the nonlinear model (M1)
and $u^h=(p^h,q^h)^T$ the approximate solution of the semilinear model (M2).
Then, the difference between the functional~$M(u)$ and~$M(u^h)$ can be approximated using
Taylor expansion, i.e., $M(u)-M(u^h)\approx M_p(u^h)(p-p^h)+M_q(u^h)(q-q^h)$. Eventually, the
first derivatives $M_p$ and $M_q$ are replaced by using the adjoint system~\eqref{adjoint_equations},
which yields after a few calculations
\begin{align}
 M(u) - M(u^h) &\; \approx\sum_{j\in \mathcal{J}_p}
 \left( \eta_{m,j}^{\text{LIN-NL}} + \eta_{t,j}^{\text{LIN}} + \eta_{x,j}^{\text{LIN}} \right)
\end{align}
with the a posteriori error estimators
\begin{subequations}
\label{eq:error_estimators}
\begin{align}
 \eta_{t,j}^{\text{LIN}} & =
   \int_0^T \! \int_{x_j^a}^{x_j^b} - \xi^T
   \big(u_t^h - R_t(u^h)\big) \, \text{d}x \, \text{d}t, \\
 \eta_{x,j}^{\text{LIN}} & =
   \int_0^T \! \int_{x_j^a}^{x_j^b} - \xi^T
   \Big(\mat{A}u^h_x - R_x(\mat{A}u^h)
   - g(u^h) + R(g(u^h))\Big) \, \text{d}x \, \text{d}t, \\
 \eta_{m,j}^{\text{LIN-NL}} & = \int_0^T \! \int_{x_j^a}^{x_j^b}
   - \xi^T B(u^h)_x \text{d}x \, \text{d}t\,.
\end{align}
\end{subequations}
Here, for the reconstruction operator $R_t$ for the temporal derivative,
we use central differences of order two, whereas $R$ and $R_x$ are defined by
\begin{subequations}
\begin{align}
R(u_i^h) =&\;\frac{1}{16}\left( -u^h_{i-3/2}+9u^h_{i-1/2}+9u^h_{i+1/2}-u^h_{i+3/2}\right),\\[2mm]
R_x(u_i^h) =&\;\frac{1}{24}\left( u^h_{i-3/2}-27u^h_{i-1/2}+27u^h_{i+1/2}-u^h_{i+3/2}\right),\;
\end{align}
\end{subequations}
which gives fourth-order accuracy with function values $u^h_{i+j/2}$ at cell centers.
The model error estimator $\eta_{m,j}^{\text{ALG-LIN}}$
between the algebraic and the semilinear model
and the discretisation error estimators
$\eta_{t,j}^{\text{NL}}$ and $\eta_{x,j}^{\text{NL}}$
for the nonlinear model
are derived analogously
for every pipe $j \in \mathcal{J}_p$;
see~\cite{Domschke2011,DomschkeKolbLang2015}.

Based on this information, the following strategy to adapt the applied models as well as the discretization parameters is used~\cite{DomschkeDuaStolwijkEtAl2018}: First, the time interval $[0,T]$ is split into certain time blocks $[T_{k-1},T_k]$ and the output functional~\eqref{equ:functional-M-general-form}, say $M_k$, and the error estimators \eqref{eq:error_estimators}, say $\eta^k_{m,j},\eta^k_{t,j},\eta^k_{x,j}$, are locally evaluated for the algebraic model and coarse meshes initially chosen. Second, models and discretization meshes are successively refined within the loop \eqref{semr_loop} until a user-prescribed tolerance $\texttt{TOL}$ is achieved in the sense that
\begin{align}
\frac{|M_k(u)-M_k(u^h)|}{|M_k(u^h)|}\approx
\frac{|\sum_{j\in \mathcal{J}_p} \left( \eta^k_{m,j} + \eta^k_{t,j} + \eta^k_{x,j} \right)|}{|M_k(u^h)|}
\le \texttt{TOL}\,.
\end{align}
The problem of finding an optimal refinement strategy is a generalisation
of the knapsack problem \cite{DomschkeDuaStolwijkEtAl2018}. We apply a greedy-like strategy
\textit{maximum error-to-cost refinement} \cite[Algorithm 3a]{DomschkeDuaStolwijkEtAl2018}
to keep the relative error below \texttt{TOL}, while retaining low computational costs.
Third, once the solution meets the tolerance requirement in $[T_{k-1},T_k]$, the
models and discretisations are coarsened if appropriate, the simulation progresses to the next
time interval, and the cycle repeats. The possible changes of the spatial discretizations between
two time blocks is treated by a conservative projection, see \cite[Sect. 6]{DomschkeKolbLang2015}.
For further details on the computation of the error estimators and the adaptive strategy, we refer to~\cite{DomschkeDuaStolwijkEtAl2018,DomschkeKolbLang2011b,DomschkeKolbLang2011a,DomschkeKolbLang2015}.

\subsection{Gradient-Based Optimization Methods}
The approach of a fully adaptive network simulation can be used to solve optimal control problems, where we consider objective functions of the generic form~\eqref{equ:functional-M-general-form}. We apply gradient-based optimization techniques, in particular the solver {\sc Donlp2}~\cite{Spellucci1998,Spellucci1998a}, where a sequential quadratic programming with an active set strategy and only equality constrained subproblems is implemented.
Since the computation of gradient information via difference quotients is rather inefficient, we apply a similar adjoint calculus as for the error estimators to get the necessary gradient information with respect to control variables of the respective optimization problem.
The whole approach has been implemented in our in-house software package {\sc Anaconda}~\cite{Kolb2011}.

Let $c$ be the vector of control variables defined by $(H_{ad}(t),\triangle p(t))$ at certain time points
and $E(u^h,c)=0$ the system of all nonlinear discretized model equations including initial, boundary and
coupling conditions. Then,
our goal is to solve the constrained optimal control problem
\begin{subequations}
\begin{align}
\text{Minimize}_{c\in C_{ad}} &\;M(u^h,c) \\
\text{subject to } &\; \quad E(u^h,c) = 0,\quad\text{where } h \text{ is chosen such that} \\
&\;\left|\sum_{j\in \mathcal{J}_p} \left( \eta^k_{m,j} + \eta^k_{t,j} + \eta^k_{x,j} \right)\right| \le |M_k(u^h,c)|\cdot \texttt{TOL}\quad\text{ for all } [T_{k-1},T_k],
\end{align}
\end{subequations}
where $C_{ad}$ is the closed and convex set of admissible controls $c$ and $M(u^h,c)$ is the output functional
defined in \eqref{equ:functional-M-general-form}, which can now also depend on the control. We would like to solve this problem with gradient-based optimization techniques. To get the necessary gradient information, we apply the same adjoint approach as described above for the error estimators. First, given a control $c\in C_{ad}$, we adaptively solve the model equations by a greedy-like refinement strategy such that the tolerance requirements are fulfilled with appropriate discretization parameters collected in $h$. We fix the discretization and solve the linear adjoint equations
\begin{align}
\left(\frac{\partial}{\partial u^h} E(u^h,c)\right)^T \xi^h =&\;-\left( \frac{\partial}{\partial u^h}M(u^h,c)\right)^T
\end{align}
for the adjoint variables $\xi^h$.
Then, the total derivative of the objective $M(u^h,c)$ with respect to the control variables $c$ reads
\begin{align}
\frac{d}{dc}M(u^h,c) =&\;\frac{\partial}{\partial c}M(u^h,c) + \left(\xi^h\right)^T
\frac{\partial}{\partial c} E(u^h,c),
\end{align}
where we have used that $\partial_{u^h}E(u^h,c)\,\partial_cu^h=-\partial_cE(u^h,c)$.
This gradient is passed to the solver {\sc Donlp2}, which delivers an improved control vector. The iteration
is stopped if a certain tolerance prescribed for {\sc Donlp2} is reached.

\section{Adaptive Simulation of Large Gas Networks} \label{sec:adaptiveErrorControl}
As an example from real gas networks in Germany, we consider the network \emph{GasLib-582} from \url{gaslib.zib.de} \cite{SchmidtAsmannBurlacuEtAl2017}. The network consists of 582 nodes (31 sources, 129 sinks, and 422 inner nodes). The nodes are connected through 278 pipes, 269 short pipes, 26 valves, 23 control valves, and 5 compressor stations, totalling 609 edges. Figure~\ref{fig:GL582-Topology} shows the network where some selected nodes are indicated that we will refer to later on.

\def\scalefig{1.15}
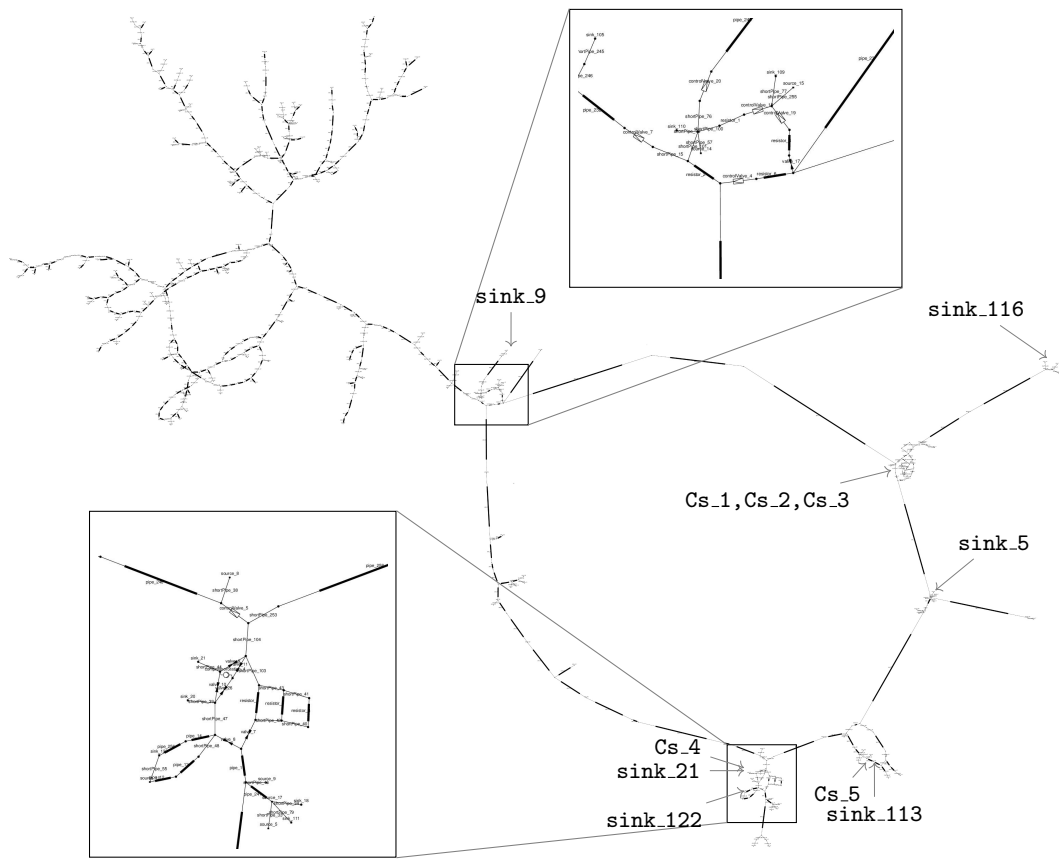
\begin{figure}[h!]
\begin{tikzpicture}[scale=\scalefig]
 \newlength{\netwidth}\setlength{\netwidth}{12.12537cm}
 \newlength{\subnetwidthOne}\setlength{\subnetwidthOne}{3.637611cm}
 \newlength{\subnetwidthTwo}\setlength{\subnetwidthTwo}{95pt}
 \node[anchor=south west] at (0,0) {\includegraphics[width=\scalefig\netwidth]{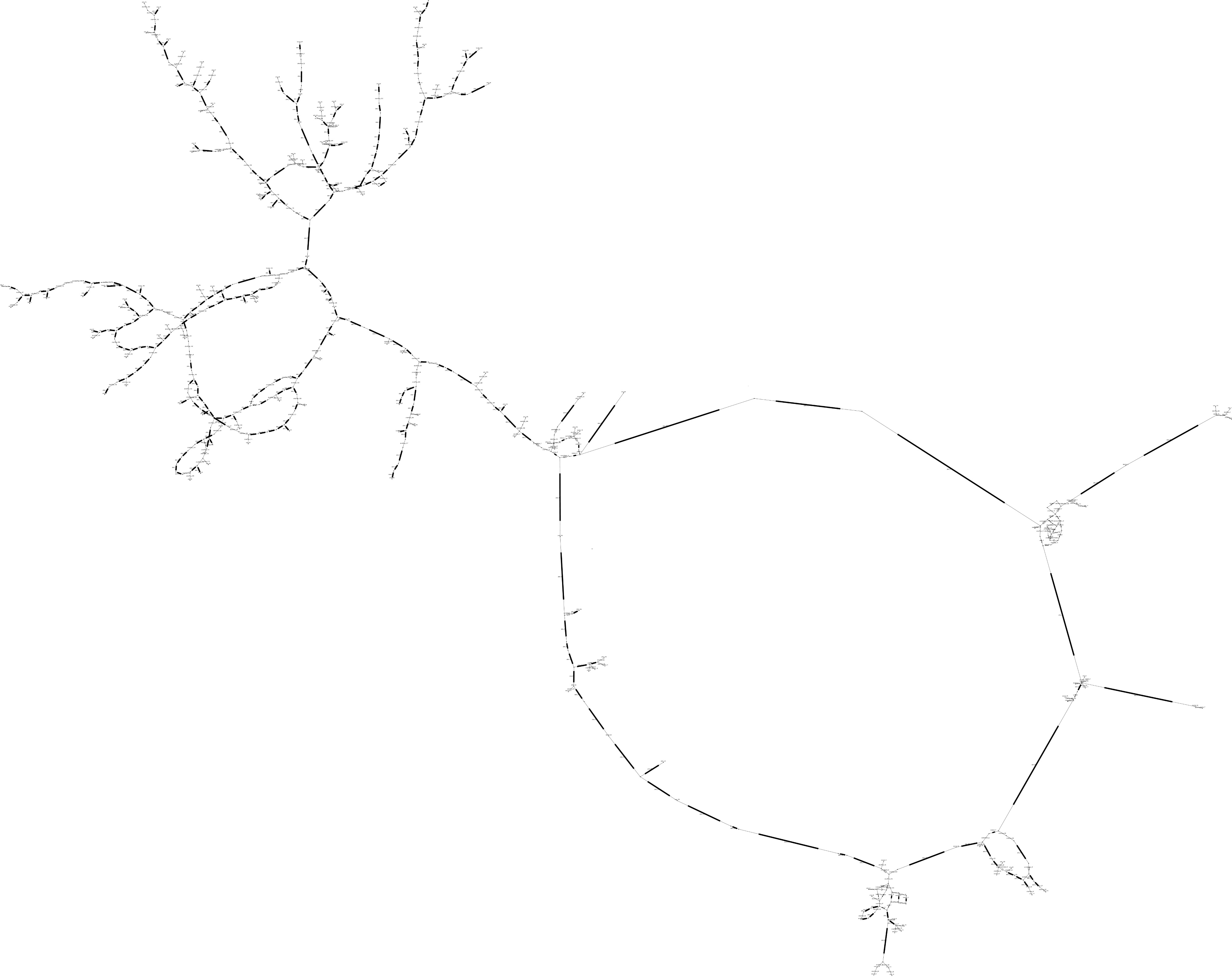}};
 \node[draw,fill=white] (A1) at (8.5,8.2) {\includegraphics[width=\scalefig\subnetwidthOne]{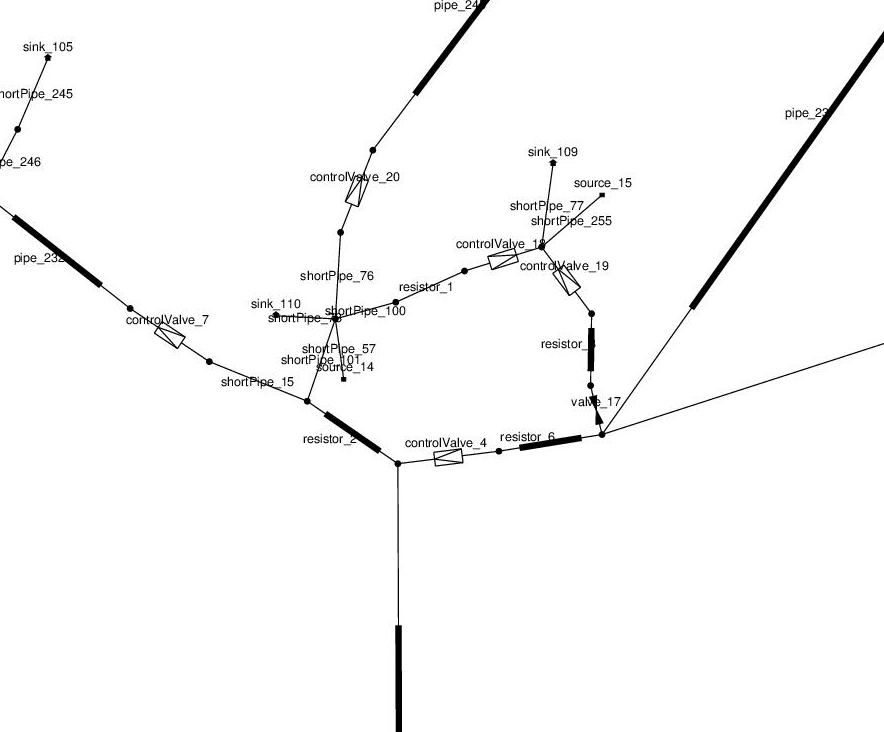}};
 \draw (5.25,5.0) rectangle (6.1,5.7);
 \draw[color=black!50] (5.25,5.7) -- (A1.north west);
 \draw[color=black!50] (6.1,5.0) -- (A1.south east);
 \node[draw,fill=white] (A3) at (2.8,2.0) {\includegraphics[width=\scalefig\subnetwidthTwo]{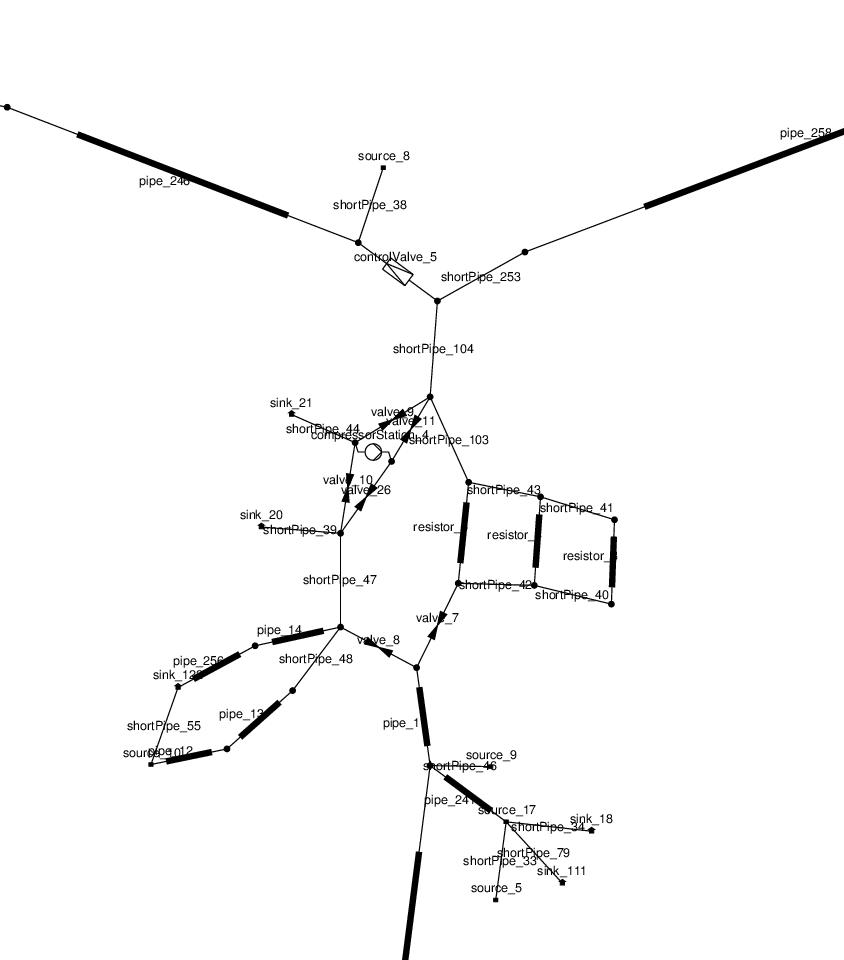}};
 \draw (8.4,0.4) rectangle (9.2,1.3);
 \draw[color=black!50] (8.4,1.3) -- (A3.north east);
 \draw[color=black!50] (8.4,0.4) -- (A3.south east);
 \node[pin={[pin edge={<-}]60:\small{\texttt{sink\_5}}}] at (10.7,3) {};
 \node[pin={[pin edge={<-}]90:\small{\texttt{sink\_9}}}] at (5.9,5.8) {};
 \node[pin={[pin edge={<-}]180:\small{\texttt{sink\_21}}}] at (8.7,1) {};
 \node[pin={[pin edge={<-}]270:\small{\texttt{sink\_113}}}] at (10.1,1.2) {};
 \node[pin={[pin edge={<-}]120:\small{\texttt{sink\_116}}}] at (12.2,5.7) {};
 \node[pin={[pin edge={<-}]200:\small{\texttt{sink\_122}}}] at (8.7,0.8) {};
 \node[pin={[pin edge={<-}]200:\small{\texttt{Cs\_1,Cs\_2,Cs\_3}}}] at (10.4,4.5) {};
 \node[pin={[pin edge={<-},pin distance=2mm]260:\small{\texttt{Cs\_5}}}] at (10.1,1.2) {};
 \node[pin={[pin edge={<-}]170:\small{\texttt{Cs\_4}}}] at (8.7,1) {};
\end{tikzpicture}
 \caption{Topology of the large real-life gas network GasLib-582.}
 \label{fig:GL582-Topology}
\end{figure}

\begin{figure}[ht]
    \centering
    \begin{subfigure}[t]{0.31\textwidth}
        \centering
        \includegraphics[width=\textwidth]{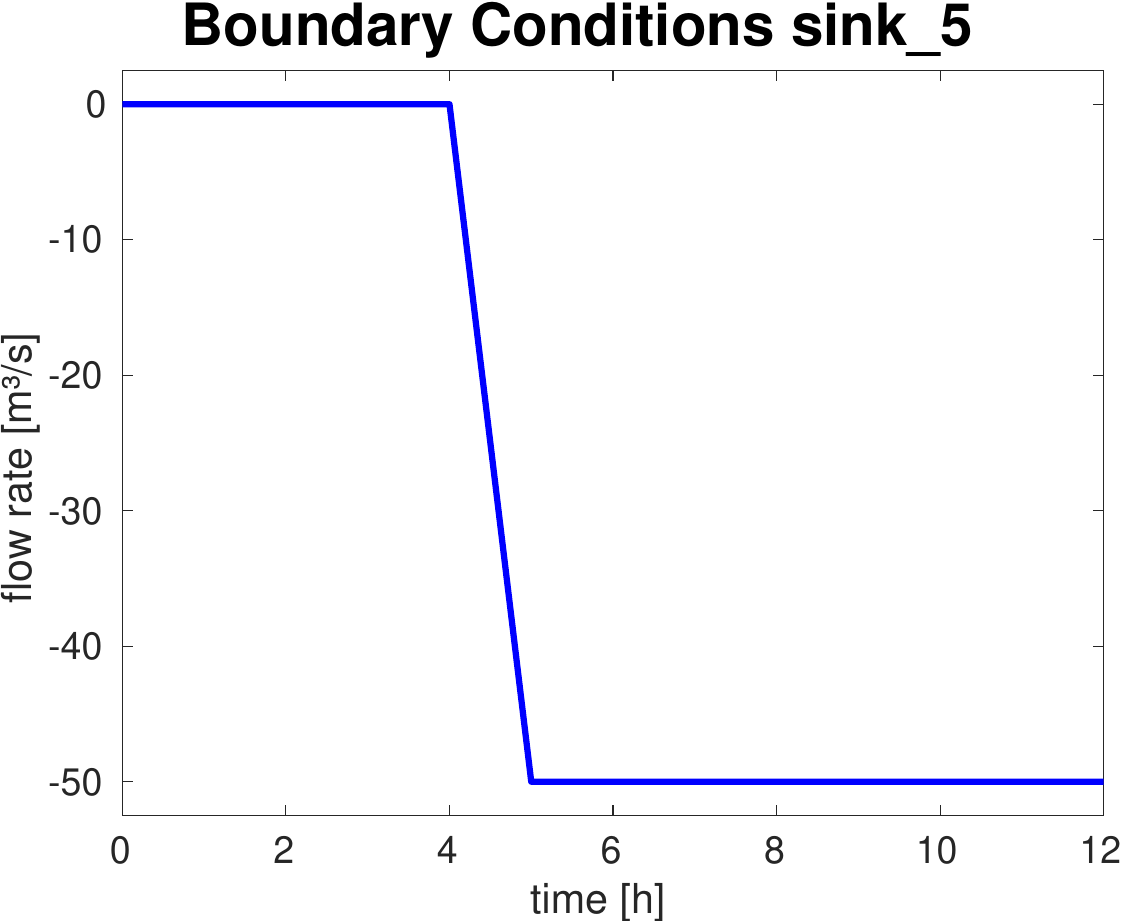}
    \end{subfigure}%
    ~
    \begin{subfigure}[t]{0.31\textwidth}
        \centering
        \includegraphics[width=\textwidth]{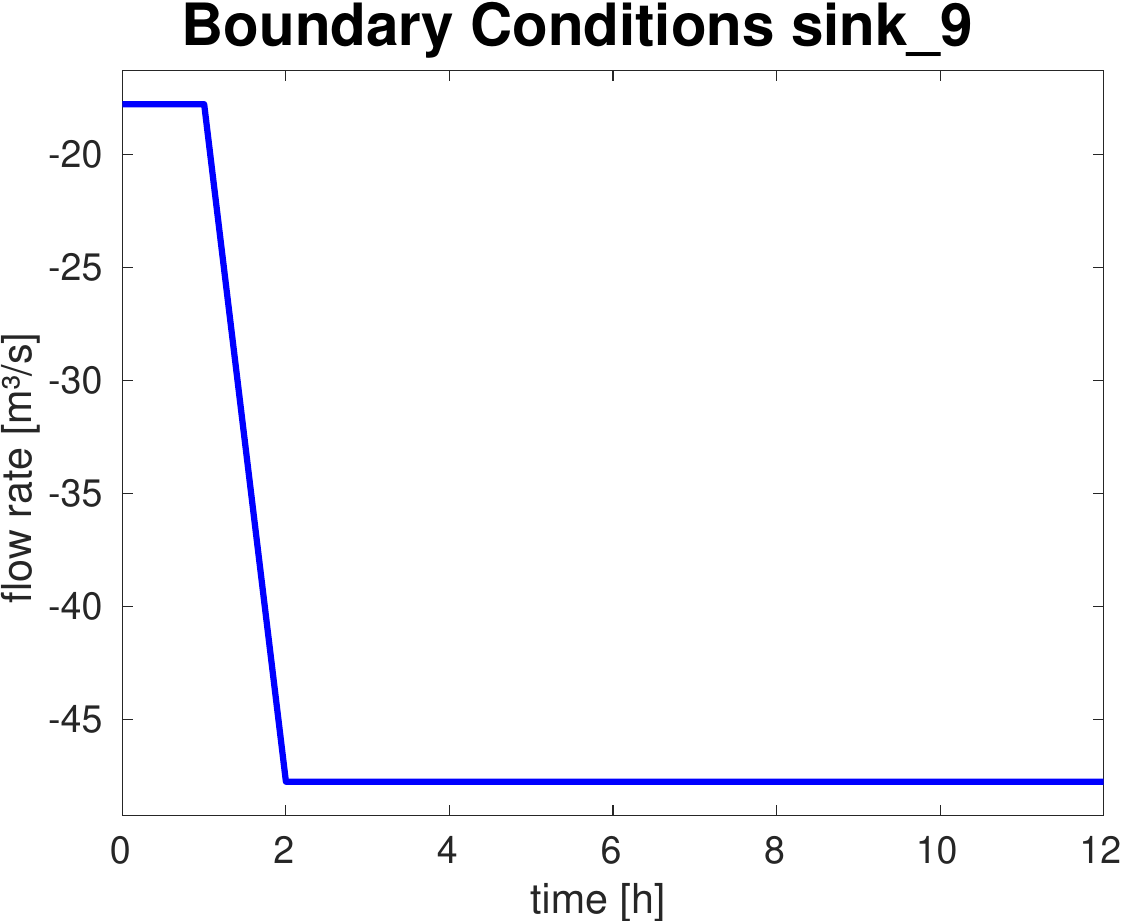}
    \end{subfigure}
    ~
    \begin{subfigure}[t]{0.31\textwidth}
        \centering
        \includegraphics[width=\textwidth]{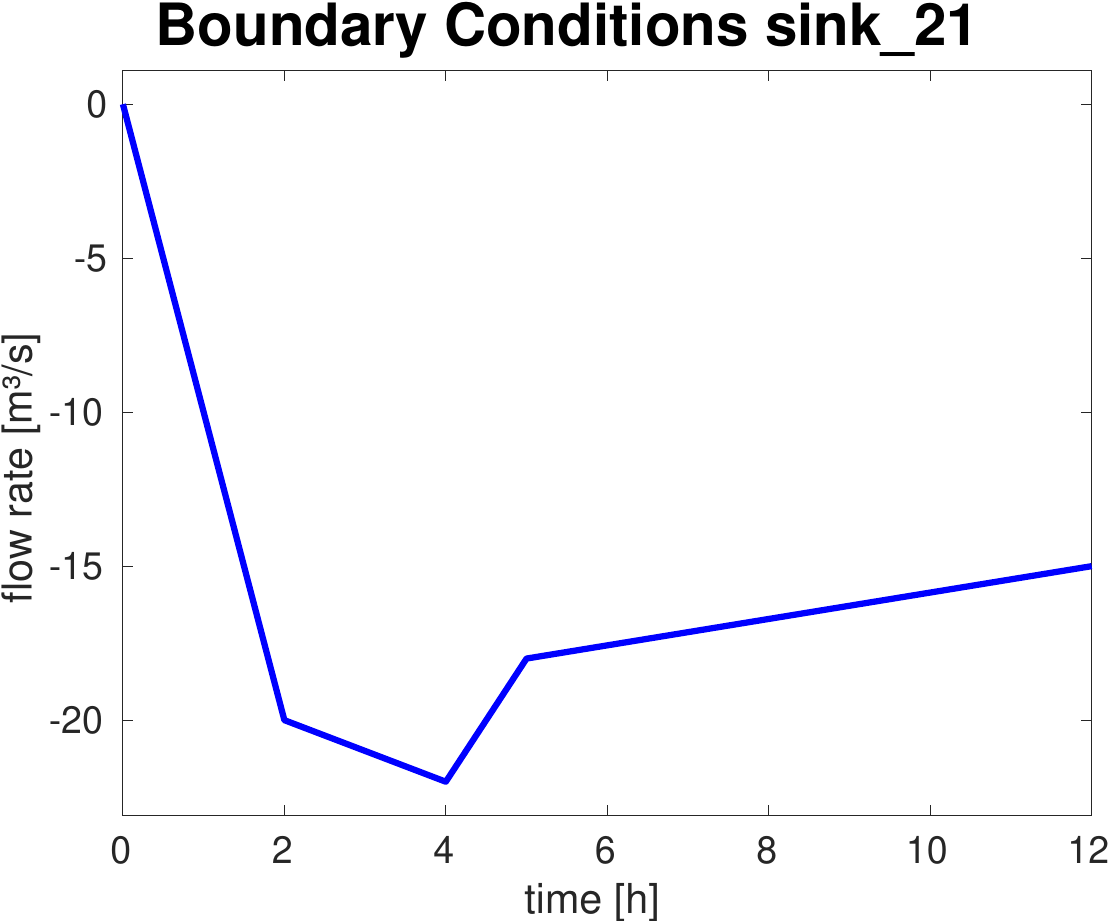}
    \end{subfigure}
    \\[7mm]
    \begin{subfigure}[t]{0.31\textwidth}
        \centering
        \includegraphics[width=\textwidth]{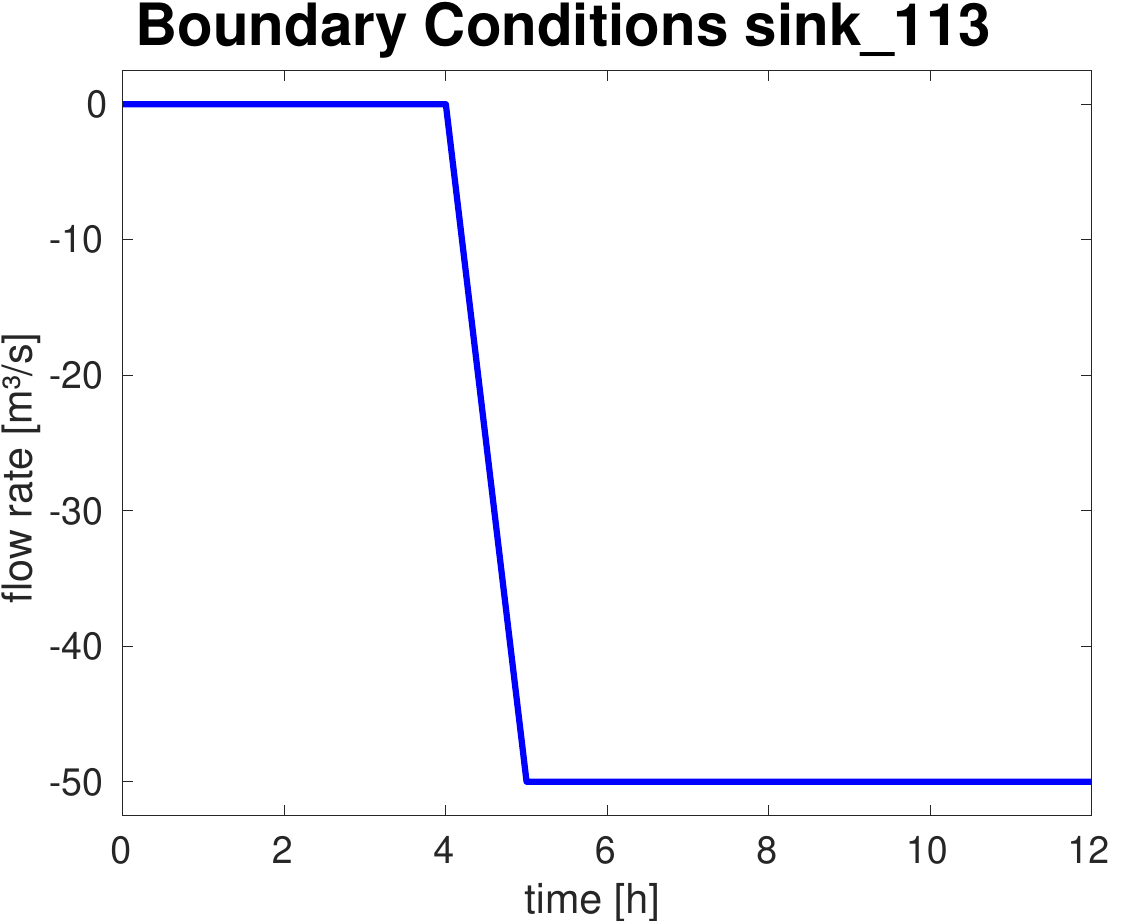}
    \end{subfigure}%
    ~
    \begin{subfigure}[t]{0.31\textwidth}
        \centering
        \includegraphics[width=\textwidth]{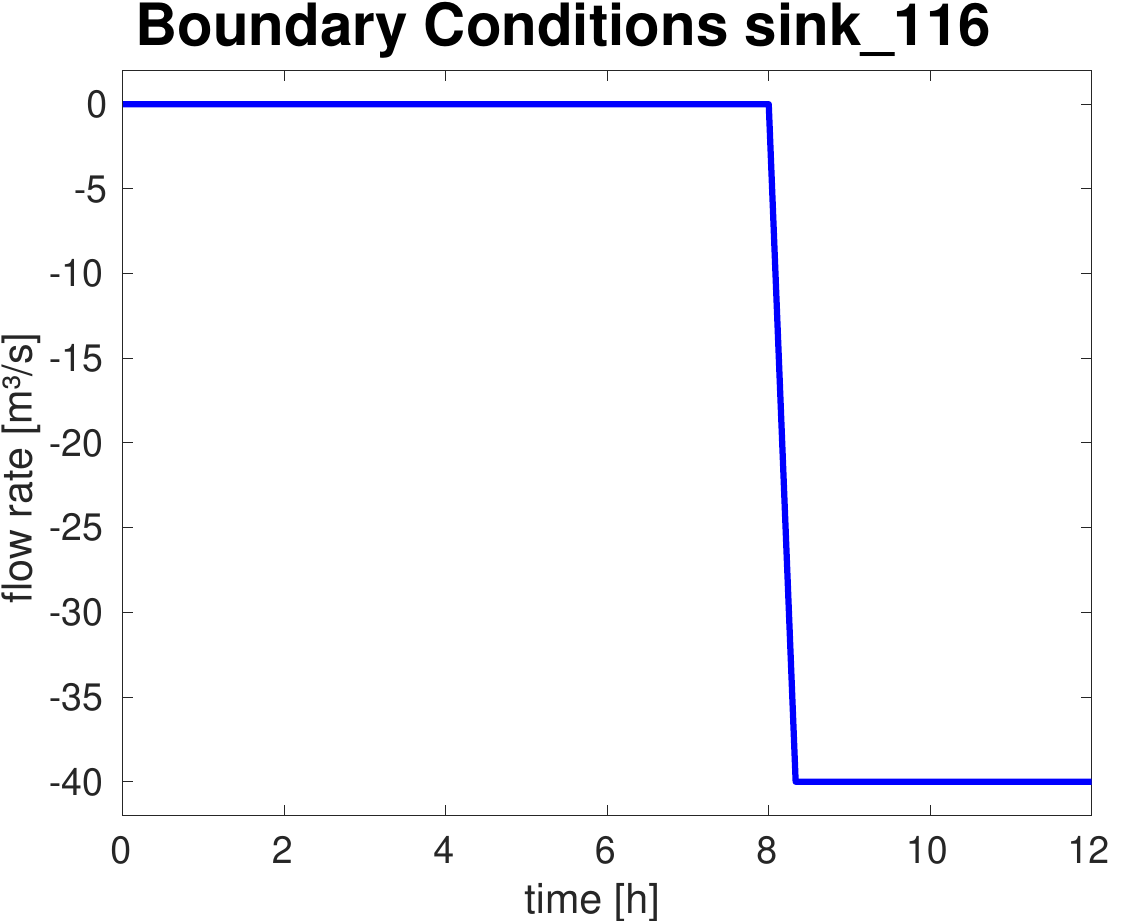}
    \end{subfigure}
    ~
    \begin{subfigure}[t]{0.31\textwidth}
        \centering
        \includegraphics[width=\textwidth]{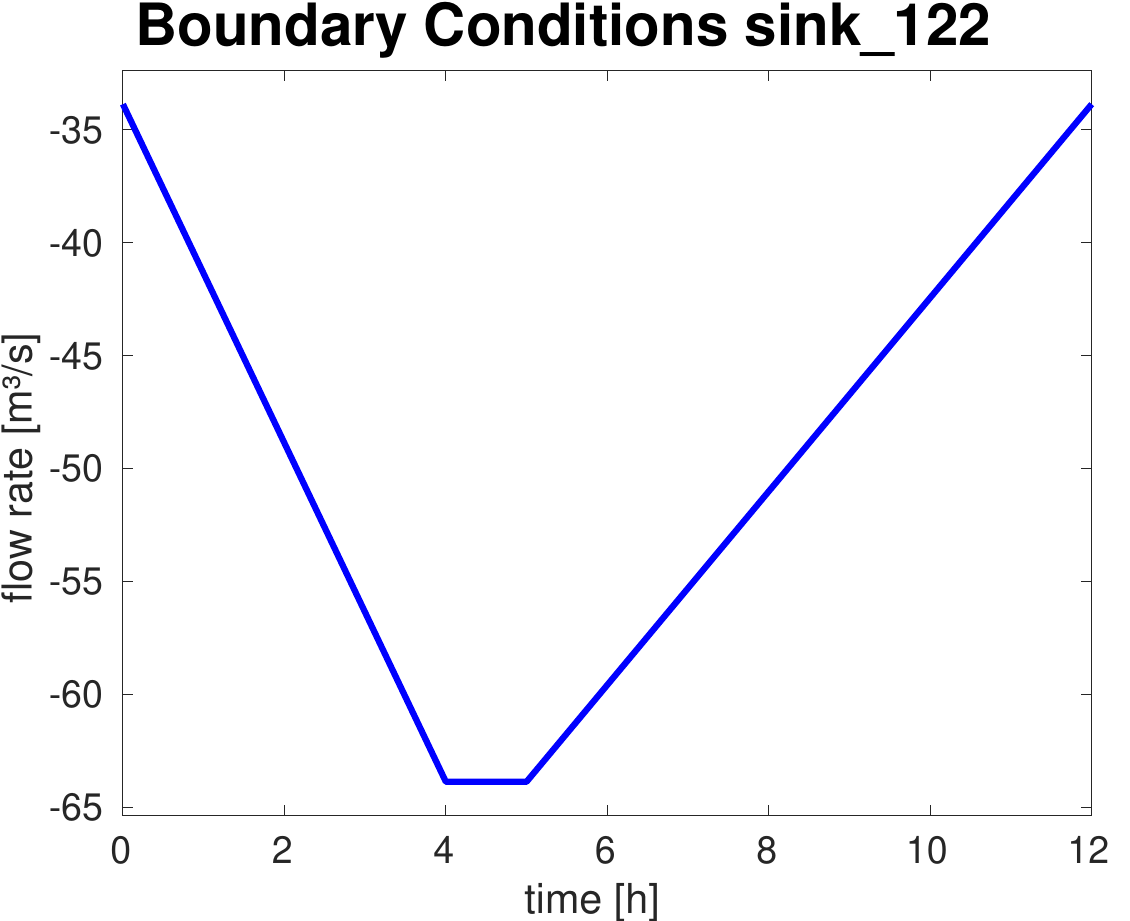}
    \end{subfigure}
    \caption{Time-dependent boundary conditions for sink 5, 9, 21, 113, 116, and 122.}
    \label{figure:GL582_BCs}
\end{figure}

\begin{figure}[htb]
    \centering
        \includegraphics[width=0.33\textwidth]{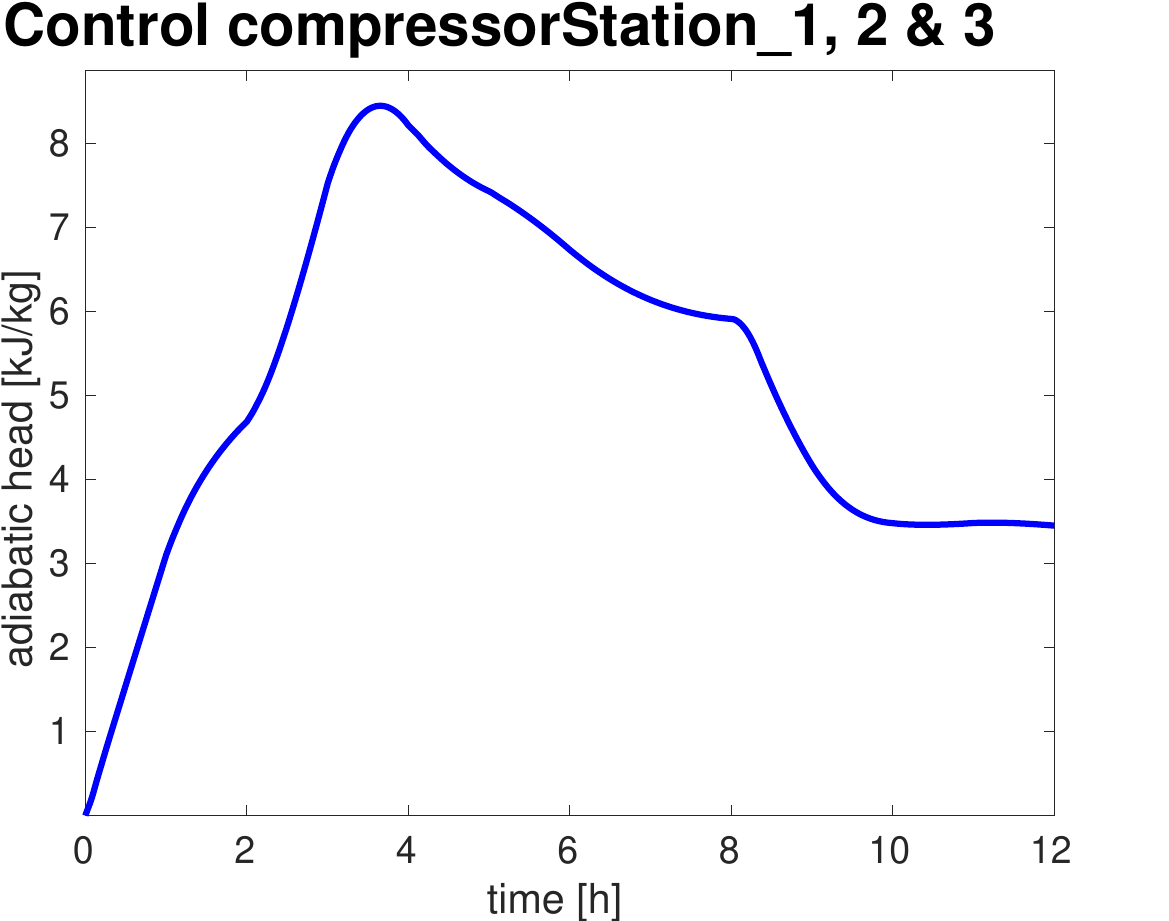}
        \includegraphics[width=0.31\textwidth]{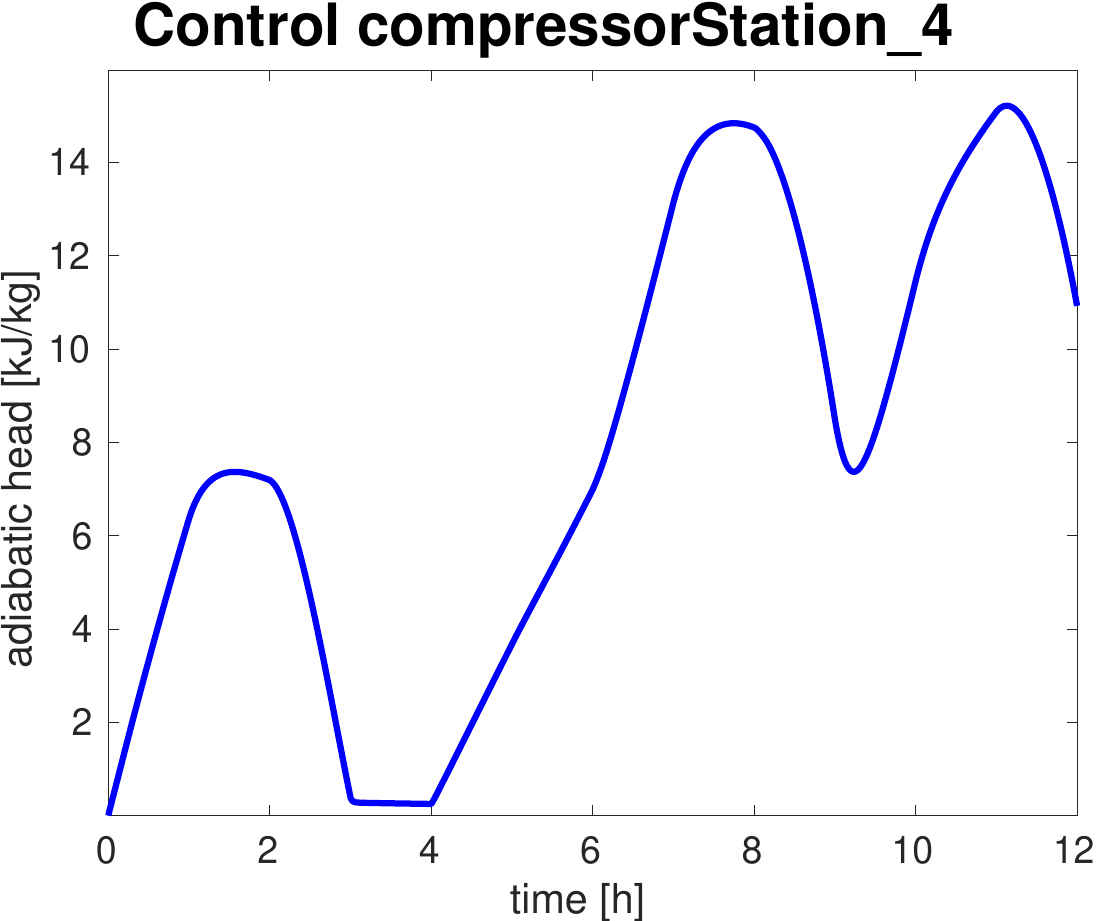}
    ~
        \includegraphics[width=0.31\textwidth]{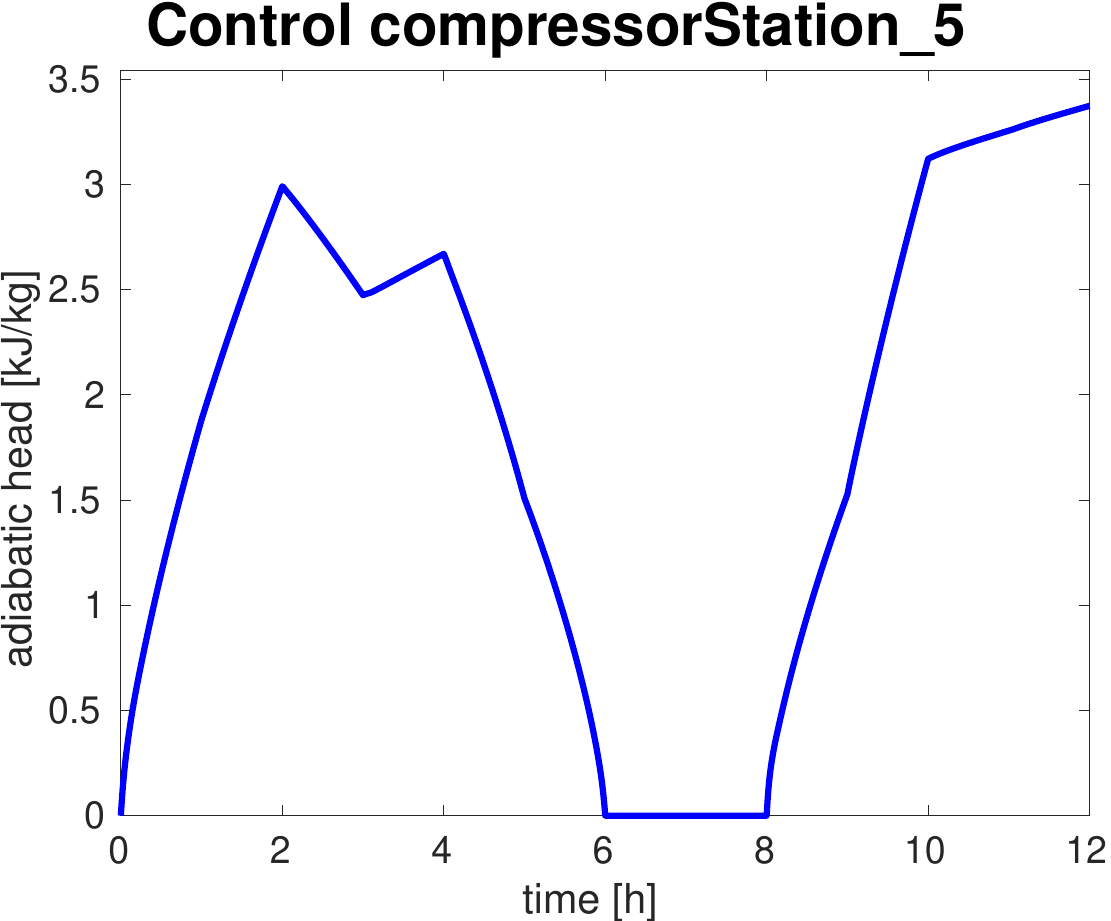}
    \caption{Time-dependent controls $\Had$ for compressors 1-5.}
    \label{figure:GL582_controlsHad}
\end{figure}

The settings of the network components with binary decisions are taken from a combined decision file originally provided by \url{gaslib.zib.de}. Initial data is generated by calculating a steady state solution for a given nomination (boundary data).
For the simulation, the boundary data of some of the sinks (\verb|sink_5|, \verb|sink_9|, \verb|sink_21|, \verb|sink_113|, \verb|sink_116|, \verb|sink_122|) are taken time-dependent as well as the controls of the compressor stations, see Figures~\ref{figure:GL582_BCs} and~\ref{figure:GL582_controlsHad}, respectively.

The adaptive simulations are run with tolerances ranging from \num{e-1} to \num{e-4}.
Running the compressor stations needs a certain amount of energy and is relatively costly. Thus, one target is to run the network at minimal compressor costs. Also, consumers like big companies nominate a certain amount of gas and may additionally have some contractual agreements relating to the pressure of the gas provided. Hence, the aim is to fulfil these requirements in order to prevent contractual penalties.
Having these examples in mind, we define the target functional to be given by the total energy consumption of the five compressor stations supplemented by the $L^2$-norm of the difference between the actual pressure and a target pressure at selected nodes.
These are $S:=\{\verb|sink_9|, \verb|sink_113|, \verb|sink_116|, \verb|sink_122|\}$, the location of which can be found in Figure~\ref{fig:GL582-Topology}. Thus, the objective function with appropriate weights reads
\begin{align}
M(u) = 0.01\cdot\sum_{cp\in\J_{cp}}\int_0^T P_{cp}(t)\;dt + \sum_{s\in S}\beta_s\,\int_0^T\left( p_s(t)-p_{s,target}(t)\right)^2\,dt\,,
\end{align}
where $\beta_s=(10^{-3},10^{-4},10^{-4},10^{-4})$ for $s=9,113,116,122$, respectively.

Table~\ref{table:GL582_SimulationResults} shows the results for the model-space-time adaptive simulations. The relative error in the target functional $M(u^h)$ with respect to a reference solution, which is computed with $\Delta t=20\,s$ and locally uniform spatial mesh sizes $\Delta x\le 250\,m$
adjusted to the pipe lengths, is presented.
We observe that the error in the target functional decreases nicely with the given tolerance. Please note that the error estimators do not give a strict upper bound, but the quality of the error estimation is quite impressive for higher tolerances \cite{Domschke2011,DomschkeKolbLang2015}.
At the same time, the CPU time needed for computation of the simulation increases moderately with the tolerance going down. Please note that for the computation of the reference solution, no adjoint equations were solved nor have any error estimators been computed. We conclude that our
fully adaptive algorithm is able to simulate a real-life gas network as GasLib-582 over a time horizon of 12 hours in a few seconds
when, e.g., a practically sufficient tolerance \num{5e-3} is requested.

\begin{table}
\centering
\caption{Simulation results for the large network GasLib-582}
\begin{tabular}{cccccc}
\hline
\TOL & rel. error & $M(u_h)$ & max/min $\Delta t$ & max/min $\Delta x$  & CPU [s] \\
\hline
5e-01 & 1.42e-02 & 837.82314 & 3600/1800 & 9936.87/255.061 & 7.06325 \\
1e-01 & 1.42e-02 & 837.82314 & 3600/1800 & 9936.87/255.061 & 8.16174 \\
5e-02 & 6.33e-03 & 831.29258 & 3600/1800 & 9936.87/255.061 & 9.46312 \\
1e-02 & 3.94e-03 & 829.32319 & 3600/1800 & 9936.87/255.061 & 11.6008 \\
5e-03 & 3.90e-03 & 829.28715 & 3600/1800 & 9936.87/255.061 & 10.8657 \\
1e-03 & 2.68e-03 & 828.28308 & 1800/450 & 9936.87/255.061 & 41.4777 \\
5e-04 & 7.54e-04 & 826.68920 & 450/450 & 9936.87/255.061 & 48.1188 \\
1e-04 & 1.43e-04 & 826.18441 & 112.5/28.125 & 9936.87/255.061 & 358.386 \\
\hline
\multicolumn{2}{c}{reference solution}  & 826.06608 & 20 & 249.984/204.049 & 533.596 \\
\hline
\end{tabular}
\label{table:GL582_SimulationResults}
\end{table}

\begin{figure}[!ht]
\begin{minipage}[c]{.49\textwidth}
\centering
\captionof{table}{Usage of models in \% for the large network}
\begin{tabular}{cccc}
\hline
\TOL & $M_1$ & $M_2$ & $M_3$ \\ \hline
1e-01 & 100.0\% & 0.0\% & 0.0\% \\
5e-02 & 95.8\% & 3.7\% & 0.5\% \\
1e-02 & 89.9\% & 7.7\% & 2.3\% \\
5e-03 & 89.1\% & 6.8\% & 4.1\% \\
1e-03 & 60.2\% & 27.6\% & 12.2\% \\
5e-04 & 41.2\% & 39.6\% & 19.1\% \\
1e-04 & 15.1\% & 65.1\% & 19.8\% \\
\hline
\end{tabular}
\label{table:GL582_models}
\end{minipage}
\begin{minipage}[c]{.49\textwidth}
  \centering
\begin{tikzpicture}[scale=0.75]
\begin{axis}[
x dir=reverse,
xmode=log,
ymode=normal,
xlabel={\texttt{TOL}},
ylabel={usage of models in \%},
ymin=0,
ymax=100,
stack plots=y,
area style,
enlarge x limits=false,
legend pos=south west
]
\addplot[color=\colALG,fill] table[x=TOL,y=model1]{modAd0_models.table} \closedcycle;
\addlegendentry{$M_1$}
\addplot[color=\colLIN,fill] table[x=TOL,y=model2]{modAd0_models.table} \closedcycle;
\addlegendentry{$M_2$}
\addplot[color=\colNL,fill] table[x=TOL,y=model3]{modAd0_models.table} \closedcycle;
\addlegendentry{$M_3$}
\end{axis}
\end{tikzpicture}
\captionof{figure}{Usage of models in \% for the simulation of the large network}
\label{figure:GL582_models}
\end{minipage}
\end{figure}

Table \ref{table:GL582_models} and Figure~\ref{figure:GL582_models} illustrate which models are used to what extent during the simulations depending on the tolerance. Not surprisingly, the smaller the tolerance, the more detailed models are used by the adaptive algorithm.

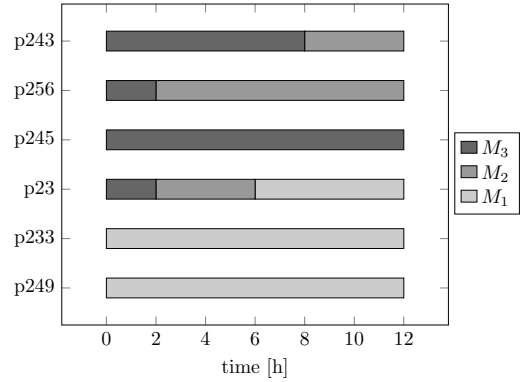
\begin{figure}[htb]
\centering
\begin{tikzpicture}[scale=0.75]
\begin{axis}[
xbar stacked,
enlargelimits=0.15,
legend style={at={(1.1,0.6)},
anchor=north,legend columns=1},
xlabel={time [h]},
symbolic y coords={p249,p233,p23,p245,p256,p243},
ytick=data,
xtick={0,2,...,12},
]
\addplot+[xbar,color=black,fill=\colNL] plot coordinates {(8,p243) (2,p256) (12,p245) (2,p23) (0,p233) (0,p249)};
\addlegendentry{$M_3$}
\addplot+[xbar,color=black,fill=\colLIN] plot coordinates {(4,p243) (10,p256) (0,p245) (4,p23) (0,p233) (0,p249)};
\addlegendentry{$M_2$}
\addplot+[xbar,color=black,fill=\colALG] plot coordinates {(0,p243) (0,p256) (0,p245) (6,p23) (12,p233) (12,p249)};
\addlegendentry{$M_1$}
\end{axis}
\end{tikzpicture}
\caption{Models used over time for a selection of pipes.\\[0.25cm]}
\label{figure:GL582_modelsOverTime}
\end{figure}

\begin{figure}[h!]
    \centering
    \begin{subfigure}[t]{0.48\textwidth}
        \centering
        \includegraphics[width=\textwidth]{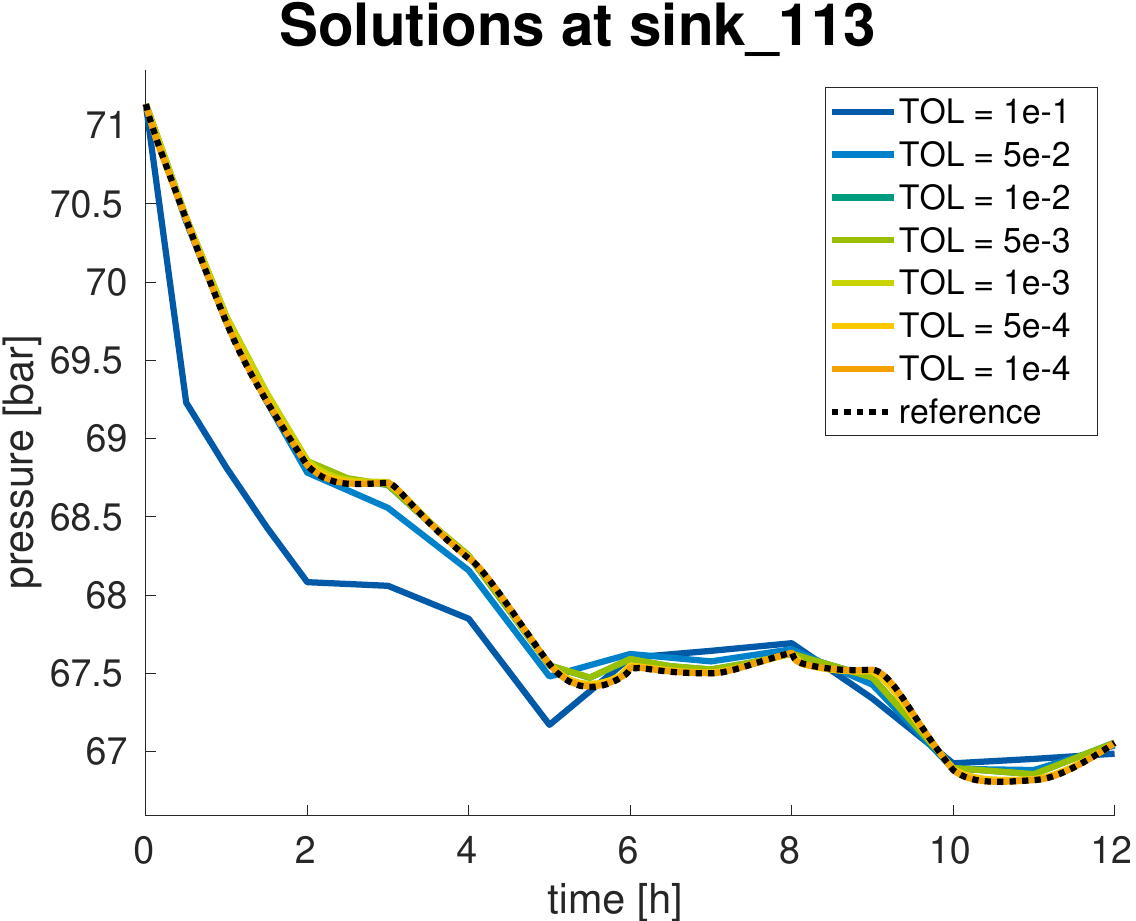}
    \end{subfigure}%
    ~
    \begin{subfigure}[t]{0.46\textwidth}
        \centering
        \includegraphics[width=\textwidth]{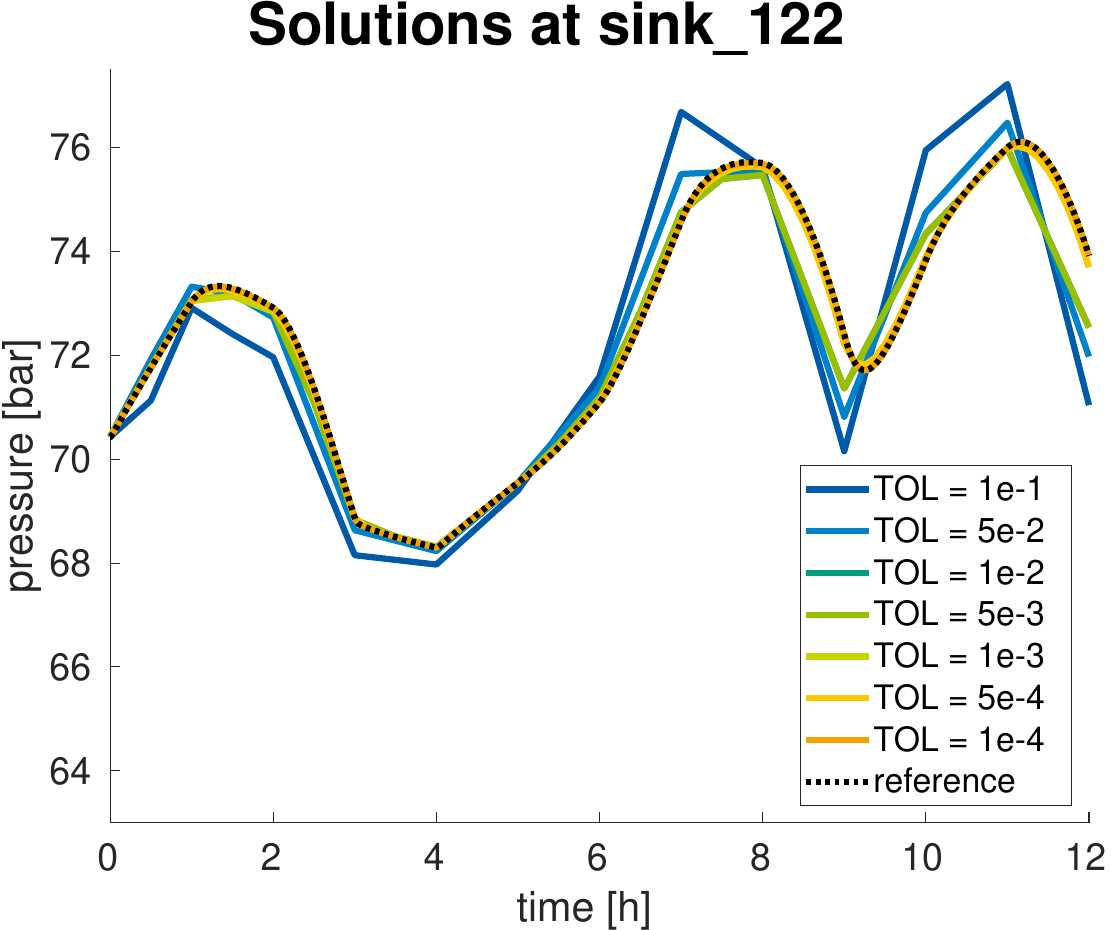}
    \end{subfigure}\\[7mm]
    \begin{subfigure}[t]{0.48\textwidth}
        \centering
        \includegraphics[width=\textwidth]{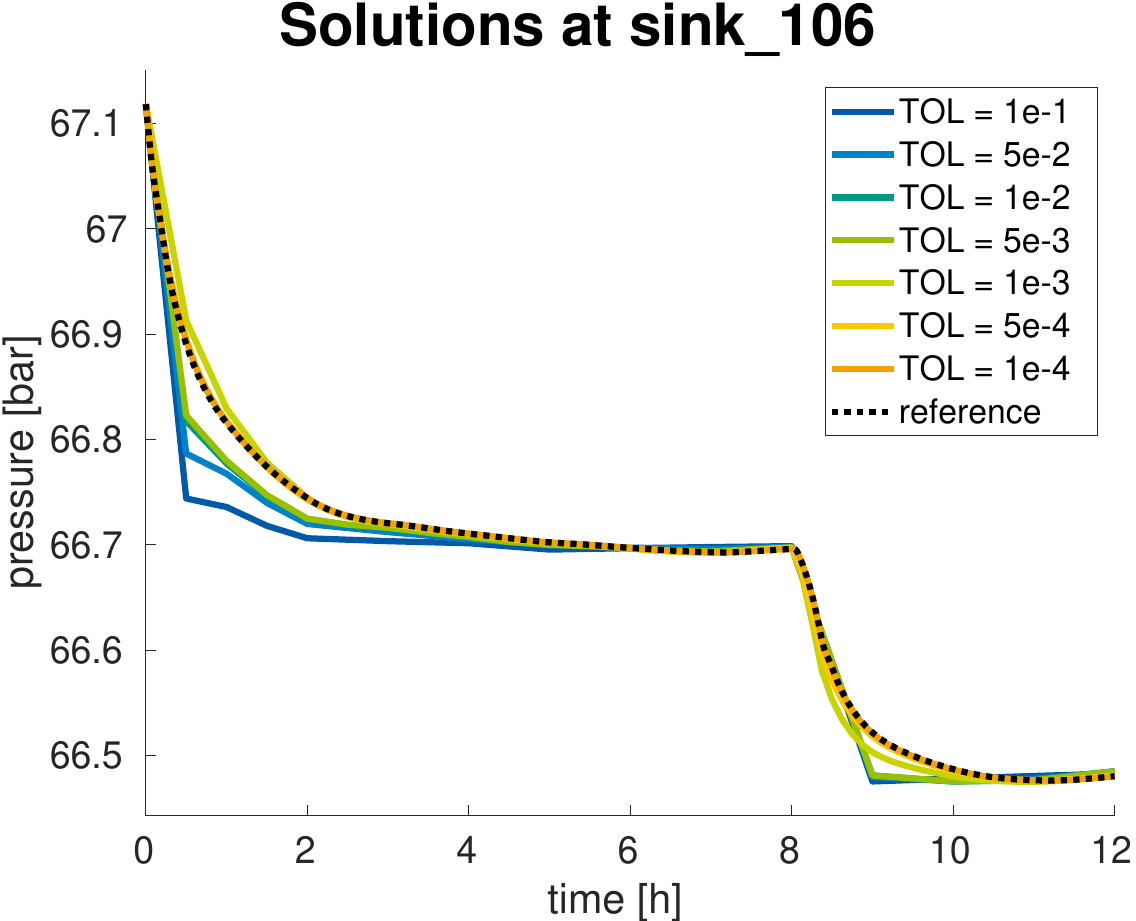}
    \end{subfigure}

    \caption{Adaptive solutions for different \TOL at sinks
    \texttt{sink\_113} and \texttt{sink\_122} contributing to the target
    functional, and \texttt{sink\_106} not contributing to the target functional.}
    \label{figure:GL582_targets}
\end{figure}

We are also interested in where and when the model refinement takes place. For a tolerance $\TOL=\num{e-3}$, we have selected a few pipes, which are adjacent to sinks in the network, see Figure~\ref{figure:GL582_modelsOverTime}. The first three pipes (p243, p256, p245) correspond to \texttt{sink\_113}, \texttt{sink\_122}, and \texttt{sink\_9}, respectively. They all contribute to the target functional. The pipe p23 corresponds to \texttt{sink\_5}, which is close to \texttt{sink\_9}. The last two pipes (p233, p249) are away from any impact to the target functional. We clearly see that pipes adjacent or close to targets use more sophisticated models than remote pipes.

Another point of interest is how well the adaptive solution corresponds to the reference solution.
Figure~\ref{figure:GL582_targets} shows the pressure at \texttt{sink\_113}, \texttt{sink\_122}, and  \texttt{sink\_106} of the adaptive simulation using different values of \TOL as well as the reference solution. It can be seen that for the sinks contributing to the target functional, i.e. \texttt{sink\_113} and \texttt{sink\_122}, the adaptive solution is quite close to the reference solution already for $\TOL=\num{e-2}$ or \num{5e-3}. For \texttt{sink\_106}, this is the case for a tolerance of \num{e-3} or even lower.

\section{Applications to Constrained Optimal Control} \label{sec:optimization}
The operation of a gas network gives rise to various scenarios which can be treated by
appropriate optimization tools. One important question in every day practice is concerned
with the issue of nomination validation:\\

\noindent
\textit{Given a stationary state $A$ of the network, is it possible
to reach a stationary state $B$ satisfying certain constraints?}\\

If the answer is yes, can we operate the network in an optimal, cost efficient way? Constraints
under consideration can be lower and upper bounds on the pressure and the flow or the operating range of a compressor station. In what follows, we exemplarily apply methods from continuous optimization to answer these questions. A key point is an appropriate modelling of compressor stations, which
is described next.

\begin{figure*}[t]
    \centering
    \begin{subfigure}[t]{0.53\textwidth}
        \centering
        \includegraphics[width=\textwidth]{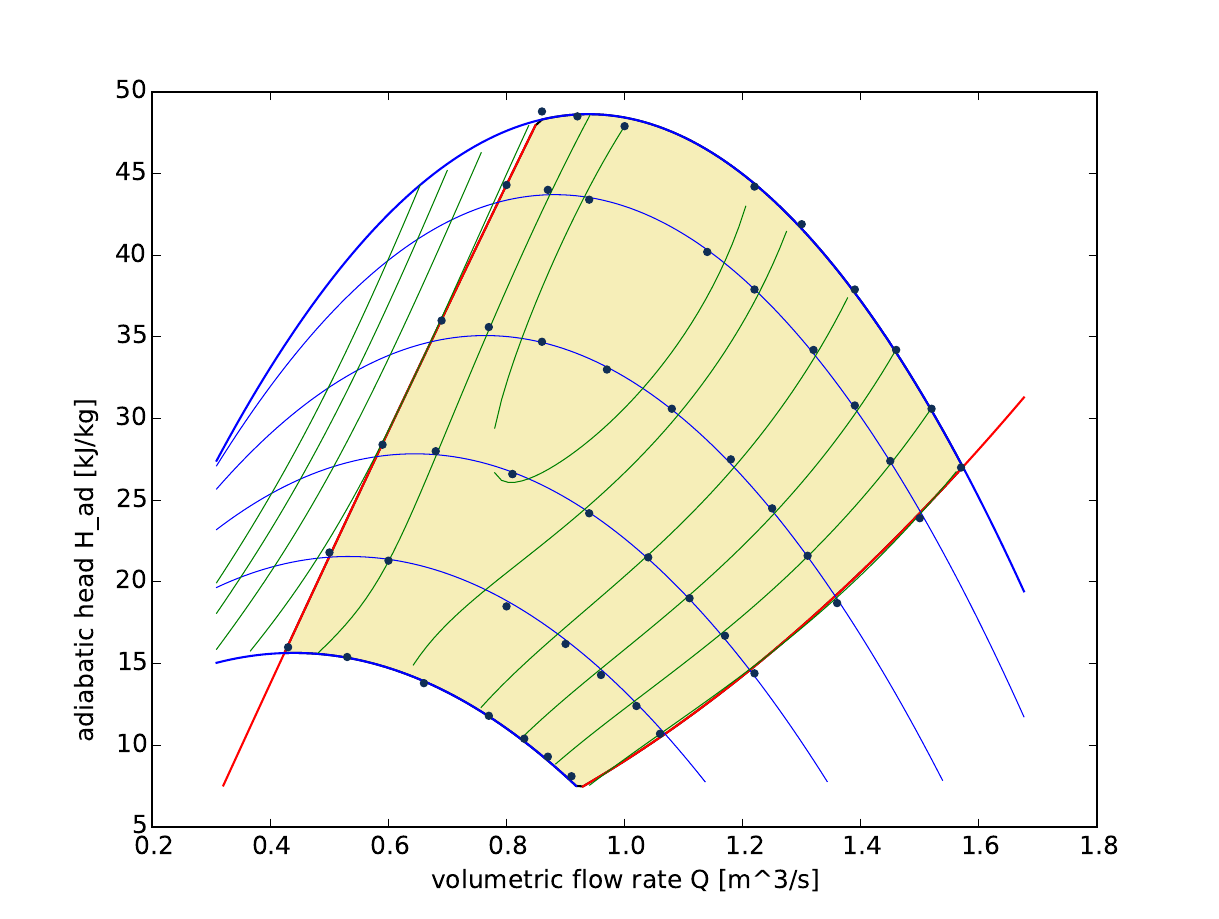}
        \caption{Single turbo compressor}\label{abb:charDiag_turbo}
    \end{subfigure}%
    ~
    \begin{subfigure}[t]{0.48\textwidth}
        \centering
        \includegraphics[width=\textwidth]{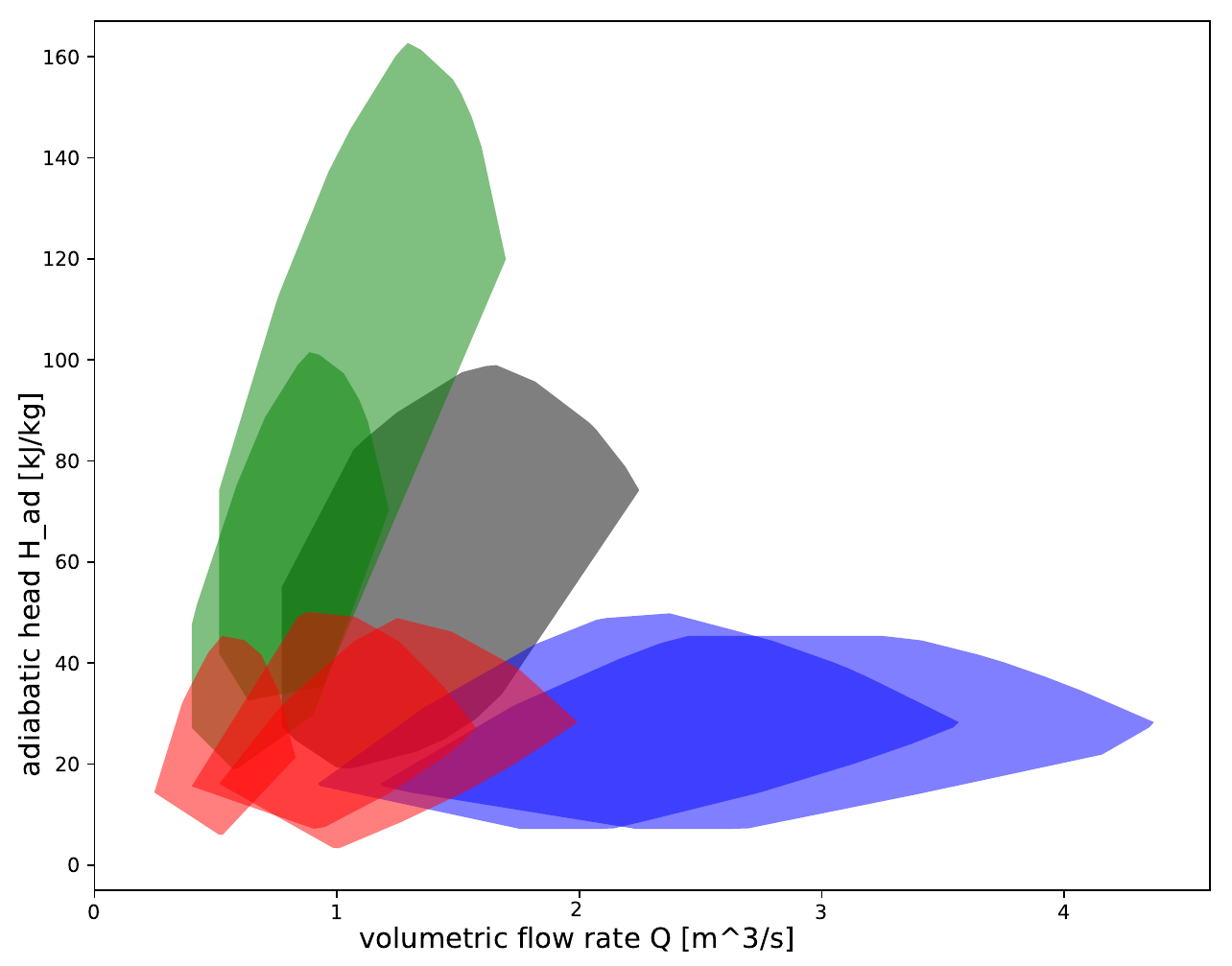}
        \caption{Aggregated characteristic field}\label{abb:charDiag_aggr}
    \end{subfigure}
    \\
    \begin{subfigure}[t]{0.48\textwidth}
        \centering
        \includegraphics[width=\textwidth]{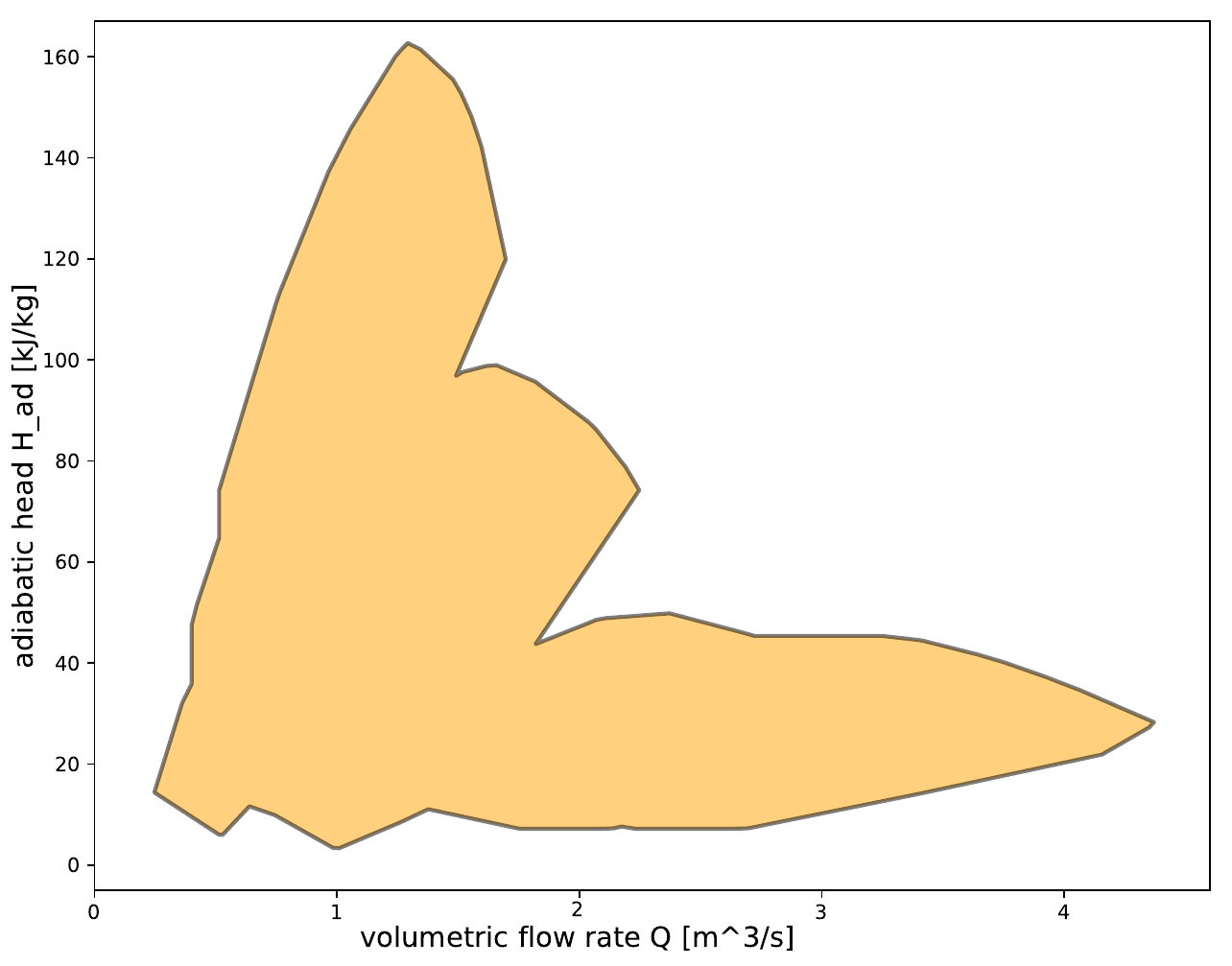}
        \caption{Approximate convex decomposition}\label{abb:charDiag_border}
    \end{subfigure}%
    ~
    \begin{subfigure}[t]{0.52\textwidth}
        \centering
        \includegraphics[width=\textwidth]{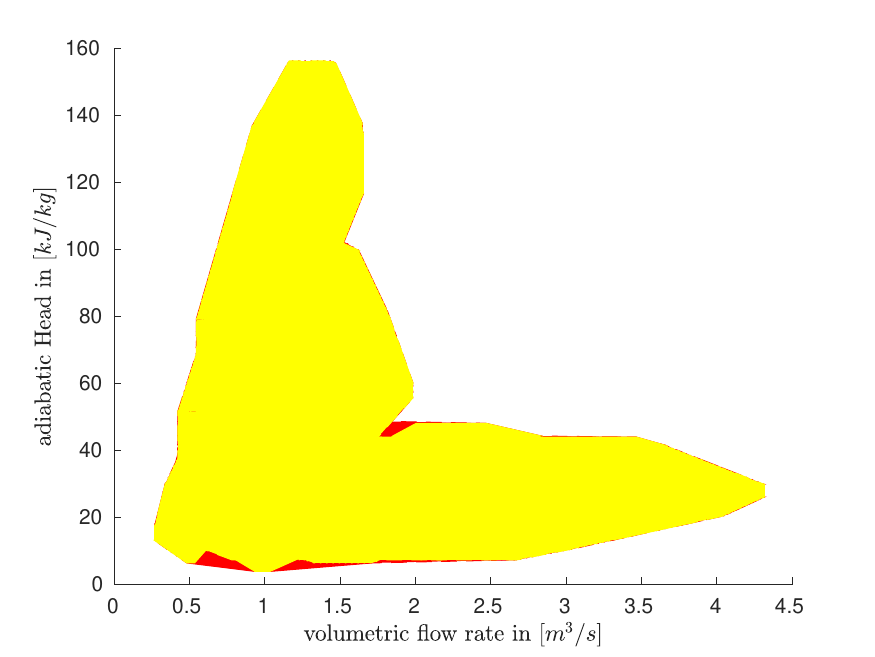}
        \caption{Semiconvex approximation}\label{abb:aggKennfeld}
    \end{subfigure}
    \caption{Modelling of compressor stations consisting of several single
    turbo compressors}
\end{figure*}

The physical model of a turbo compressor is determined by a characteristic field \cite{WaltherHiller2017}, see Figure~\ref{abb:charDiag_turbo}. A compressor station typically consists of multiple compressors that can be run in different configurations: single, serial, and parallel. Each configuration is again described by a characteristic field that overlap largely, see Figure~\ref{abb:charDiag_aggr}. Resulting from approximate convex decomposition \cite{HillerWalther2017,LienAmato2006}, we take an outer linear approximation of the physical model of the aggregated characteristic fields as a base, see Figure~\ref{abb:charDiag_border}. In order to make such a model applicable within an optimization, we approximate it by a semiconvex set $K$ in the $Q$-$\Had$-space, see Figure~\ref{abb:aggKennfeld}. The approximation is semiconvex in the sense that for fixed $\Had$, if $(Q_1,\Had)\in K$ and $(Q_2,\Had)\in K$, then
\begin{align}
\label{cond:semiconvex}
(\lambda Q_1 + (1-\lambda) Q_2,H_{ad}) \in K \qquad \forall \lambda \in [0,1]\,.
\end{align}
The difference between the approximate convex decomposition model and the semiconvex approximation is displayed in Figure~\ref{abb:aggKennfeld}. Property \eqref{cond:semiconvex} allows to incorporate the characteristic
diagram of a compressor station as a state constraint for the volumetric flow rate $Q$ through it,
see \eqref{opt:condQh} below. The pressure head $H_{ad}$ is then only restricted by a flow-independent
minimum and maximum value, $H_{ad,min}$ and $H_{ad,max}$, respectively.

The path from A to B mainly depends on the boundary values, especially the change of the flow
rate $q_s$ at the sinks. One
could prescribe this change by a certain function in time, e.g. a linear one, over an appropriate time horizon and
try to find a set of controls that guarantees the compliance of all operational restrictions. However, the identification of such
transfer functions is often a difficult task. A second option, which we will follow here, is to relax the inflexibility
of a fixed time-dependent boundary condition by considering the outflow rates $q_s$, $s\in\J_q$, as part of the controls and
adding a weighted tracking type functional to the objective function, resulting in
\begin{align}\label{obj:Mnv}
M_{nv}(u,c) = \alpha\cdot\sum_{cp\in\J_{cp}}\int_0^T P_{cp}(t)\;dt + \sum_{s\in\J_q}\beta_s\,\int_0^T\left( q_s(t)-q_{s,target}(t)\right)^2\,dt\,.
\end{align}
with positive weights $\alpha$, $\beta_s$, and $c=(\{H_{ad}^{cp}\},\{q_s\})$ - the set of all pressure heads and outflow rates that have to be changed. It is beneficial that $q_{s,target}(t)=q_{s}^B$ for $t\in [T-t^\star,T]$ for a sufficiently large $t^\star<T$ and $q_{s}^B$ being the desired new outflow rate at state $B$.
The admissible set is defined as
\begin{align}
C_{ad}:=\{(H_{ad}^{cp},q_s)_{cp\in\J_{cp},s\in\J_q}:H_{ad,min}^{cp}\le H_{ad}^{cp} \le H_{ad,max}^{cp},
\;q_s\in\R \}.
\end{align}
We also consider lower and upper bounds for the flow rate $q_s$ and the pressure $p_s$ at certain sinks
and sources with index sets $\J_{pl}$, $\J_{pu}$, $\J_{ql}$, and $\J_{qu}$, respectively. In general, $\J_{pl}\cap\J_{pu}\neq\emptyset$ and $\J_{ql}\cap\J_{qu}\neq\emptyset$ .

Our resulting optimal control problem for the nomination validation reads as follows:
\begin{subequations}
\begin{align}
\text{Minimize}_{c\in C_{ad}} &\;M_{nv}(u^h,c) \\
\text{subject to } &\; \quad E(u^h,c) = 0,\quad\text{where } h \text{ is chosen such that} \\
&\;\left|\sum_{j\in \mathcal{J}_p} \left( \eta^k_{m,j} + \eta^k_{t,j} + \eta^k_{x,j} \right)\right| \le |M_k(u^h,c)|\cdot \texttt{TOL}\quad\text{ for all } [T_{k-1},T_k], \\[2mm]
&\; Q_{min}(H_{ad}^{cp}) \le Q^h_{cp} \le Q_{max}(H_{ad}^{cp})
\quad\text{for all }cp\in\J_{cp}, \label{opt:condQh}\\[2mm]
&\; p^h_s(t_i)\ge p^l_s(t_i) \quad\text{for all discrete time points }t_i \text{ and } s\in\J_{pl},\\[1mm]
&\; p^h_s(t_i)\le p^u_s(t_i) \quad\text{for all discrete time points }t_i \text{ and } s\in\J_{pu},\label{opt:condPr}\\[1mm]
&\; q^h_s(t_i)\ge q^l_s(t_i) \quad\text{for all discrete time points }t_i \text{ and } s\in\J_{ql},\\[1mm]
&\; q^h_s(t_i)\le q^u_s(t_i) \quad\text{for all discrete time points }t_i \text{ and } s\in\J_{qu}\label{opt:condQ}\,.
\end{align}
\end{subequations}
We apply the gradient-based SQP-solver {\sc Donlp2}~\cite{Spellucci1998,Spellucci1998a},
which can handle control as well as
state constraints of the form \eqref{opt:condQh}-\eqref{opt:condQ}.

In the following sections, we show the applicability of the continuous optimization algorithm for nomination validation including the semiconvex approximation of the characteristic diagram as constraints for the compressor stations.

\subsection{Network with Three Compressor Stations}
\begin{figure}[!ht]
\centering 
\begin{tikzpicture}[y=-1cm,scale=0.09]
\draw[thick,black] (-51.62444,25.43333) -- (-41.04,9.55778);
\draw[line width=3.14pt,black] (-49.50667,22.25778) -- (-43.15778,12.73333);
\draw[thick,black] (-19.87333,9.55778) -- (-9.29111,-1.02444);
\draw[line width=3.14pt,black] (-17.75778,7.44222) -- (-11.40667,1.09111);
\draw[thick,black] (-9.29111,-1.02444) -- (1.29333,9.55778);
\draw[line width=3.14pt,black] (-7.17333,1.09111) -- (-0.82444,7.44222);
\draw[thick,black] (-19.87333,9.55778) -- (-9.29111,20.14222);
\draw[line width=3.14pt,black] (-17.75778,11.67556) -- (-11.40667,18.02444);
\draw[thick,black] (-9.29111,20.14222) -- (1.29333,9.55778);
\draw[line width=3.14pt,black] (-7.17333,18.02444) -- (-0.82444,11.67556);
\draw[thick,black] (22.46,9.55778) -- (33.04222,-1.02444);
\draw[line width=3.14pt,black] (24.57556,7.44222) -- (30.92667,1.09111);
\draw[thick,black] (33.04222,-1.02444) -- (43.62667,9.55778);
\draw[line width=3.14pt,black] (35.16,1.09111) -- (41.50889,7.44222);
\draw[thick,black] (22.46,9.55778) -- (33.04222,20.14222);
\draw[line width=3.14pt,black] (24.57556,11.67556) -- (30.92667,18.02444);
\draw[thick,black] (33.04222,20.14222) -- (43.62667,9.55778);
\draw[line width=3.14pt,black] (35.16,18.02444) -- (41.50889,11.67556);
\draw[thick,black] (54.20889,20.14222) -- (64.79333,20.14222);
\draw[line width=3.14pt,black] (56.32667,20.14222) -- (62.67556,20.14222);
\draw[thick,black] (64.79333,9.55778) -- (64.79333,-6.31556);
\draw[line width=3.14pt,black] (64.79333,6.38444) -- (64.79333,-3.14222);
\draw[thick,black] (-51.62444,-6.31556) -- (-41.04,9.55778);
\draw[thick,black] (-9.29111,-1.02444) -- (0.58667,-1.02444);
\draw[thick,black] (1.99778,-1.02444) -- (11.87556,-1.02444);
\pic at (1.29333,-1.02444) [scale=0.7] {valve};
\path (0.97556,-1.22) node[text=black,anchor=base west] {Vl};
\draw[thick,black] (33.04222,20.14222) -- (42.92,20.14222);
\pic at (43.6,20.14222) [scale=0.7] {controlValve};
\draw[thick,black] (44.33111,20.14222) -- (54.20889,20.14222);
\path (43.30889,19.59556) node[text=black,anchor=base west] {Cv};
\draw[thick,black] (-41.04,9.55778) -- (-19.87333,9.55778);
\pic at (-30.45778,9.55778) [scale=0.7] {compressor};
\path (-31.11111,8.41556) node[text=black,anchor=base west] {Cs1};
\draw[thick,black] (1.29333,9.55778) -- (22.46,9.55778);
\pic at (11.87556,9.55778) [scale=0.7] {compressor};
\path (11.22222,8.41556) node[text=black,anchor=base west] {Cs2};
\draw[thick,black] (43.62667,9.55778) -- (64.79333,9.55778);
\pic at (54.20889,9.55778) [scale=0.7] {compressor};
\path (53.55556,8.41556) node[text=black,anchor=base west] {Cs3};
\pic at (-51.72889,25.52222) [scale=0.7] {source};
\path (-51.94222,25.02667) node[text=black,anchor=base east] {S1};
\pic at (-51.72889,-6.22889) [scale=0.7] {source};
\path (-51.94222,-6.72222) node[text=black,anchor=base east] {S2};
\pic at (64.47556,19.73556) [scale=0.6,rotate=180] {sink};
\path (64.47556,19.73556) node[text=black,anchor=base west] {T3};
\pic at (64.47556,-6.72222) [scale=0.6,rotate=180] {sink};
\path (64.47556,-6.72222) node[text=black,anchor=base west] {T2};
\pic at (11.98222,-1.11333) [scale=0.6,rotate=180] {sink};
\path (11.71778,-1.43111) node[text=black,anchor=base west] {F};
\pic at (-9.43111,20.05333) [scale=0.6,rotate=180] {sink};
\path (-9.60889,19.73556) node[text=black,anchor=north west] {T1};
\filldraw[thick,black] (-41.04,9.55778) circle (0.14222cm);
\filldraw[thick,black] (-19.87333,9.55778) circle (0.14222cm);
\filldraw[thick,black] (-9.29111,-1.02444) circle (0.14222cm);
\filldraw[thick,black] (1.29333,9.55778) circle (0.14222cm);
\filldraw[thick,black] (22.46,9.55778) circle (0.14222cm);
\filldraw[thick,black] (33.04222,-1.02444) circle (0.14222cm);
\filldraw[thick,black] (33.04222,20.14222) circle (0.14222cm);
\filldraw[thick,black] (54.20889,20.14222) circle (0.14222cm);
\filldraw[thick,black] (43.62667,9.55778) circle (0.14222cm);
\filldraw[thick,black] (64.79333,9.55778) circle (0.14222cm);

\end{tikzpicture}%
\caption{Network with three compressor stations}
\label{figure:Steinach}
\end{figure}
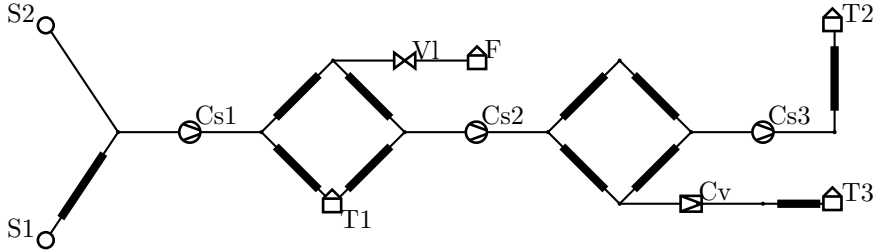

For demonstration purpose, we firstly consider a simplified, but very illustrative
network with three compressor stations, see Fig.~\ref{figure:Steinach}.

\begin{figure*}[t!]
    \centering
    \begin{subfigure}[t]{0.42\textwidth}
        \centering
        \includegraphics[width=\textwidth]{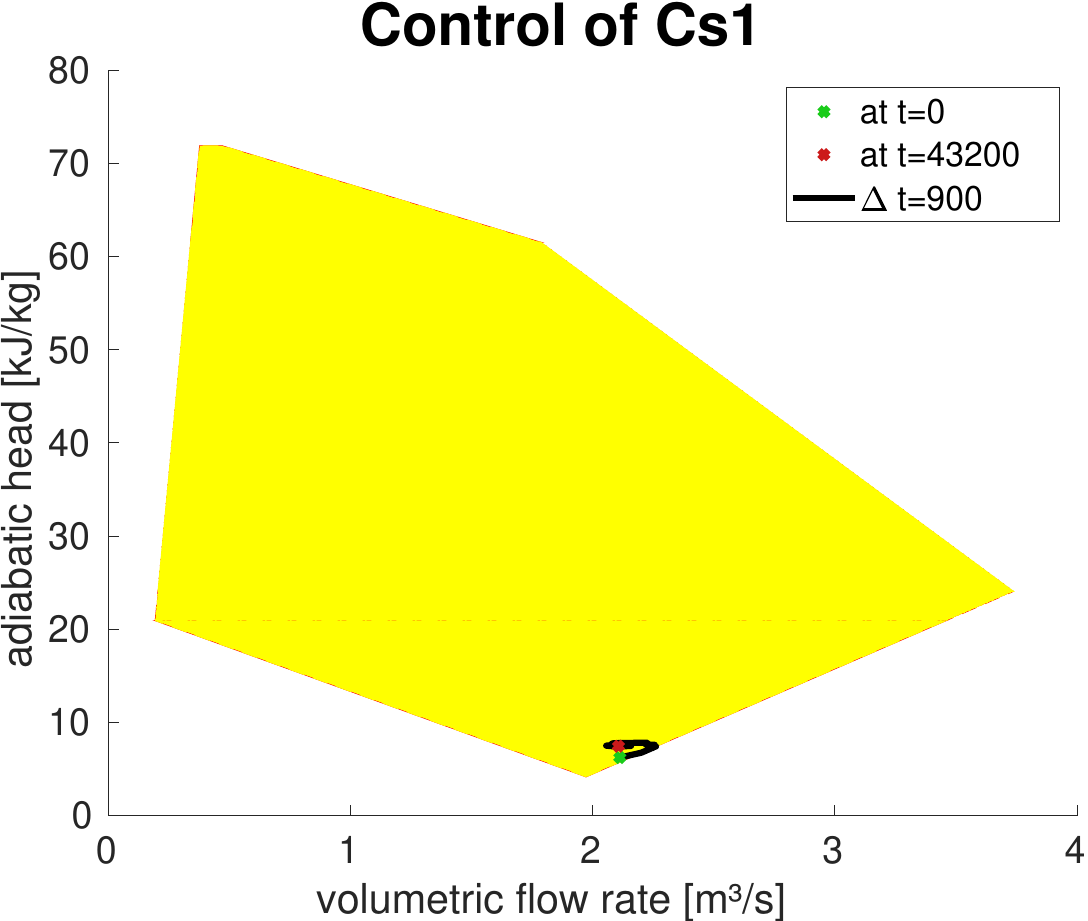}
    \end{subfigure}%
    ~
    \begin{subfigure}[t]{0.43\textwidth}
        \centering
        \includegraphics[width=\textwidth]{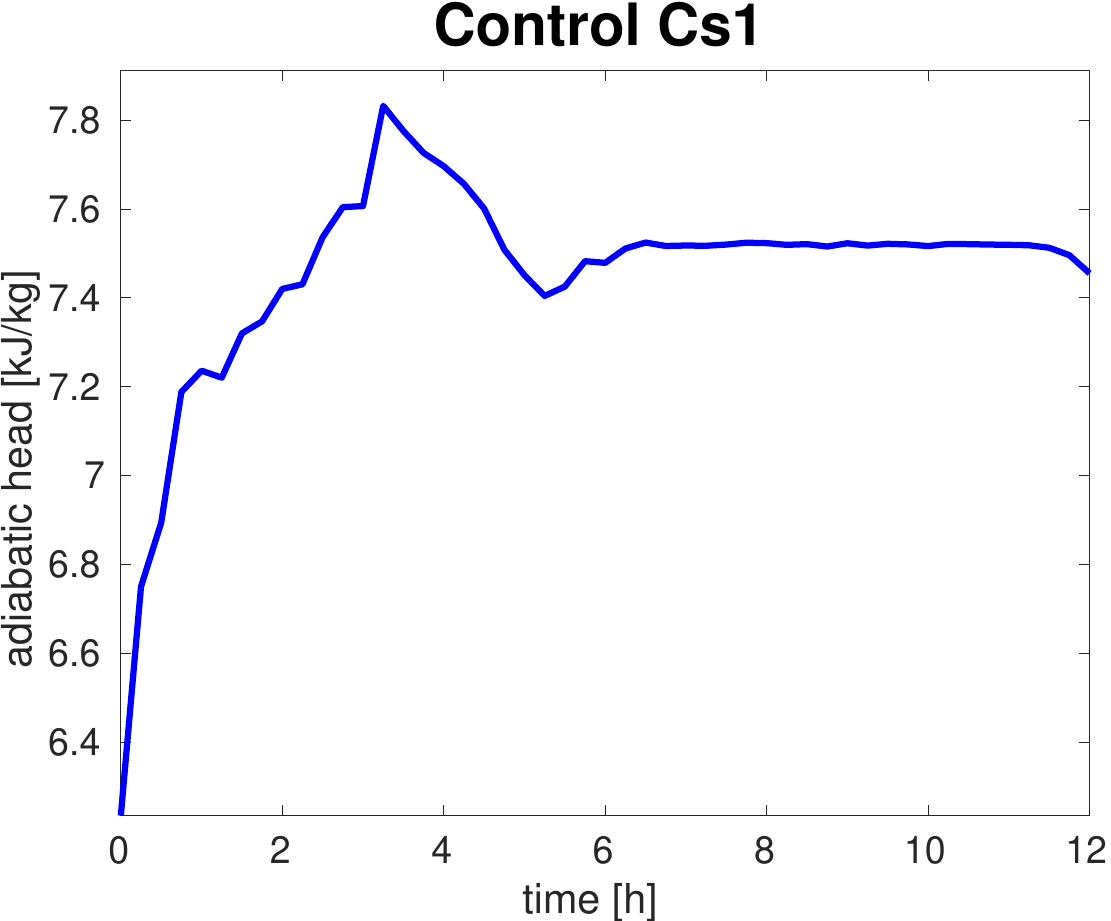}
    \end{subfigure}
    \\[0.5cm]
    \begin{subfigure}[t]{0.42\textwidth}
        \centering
        \includegraphics[width=\textwidth]{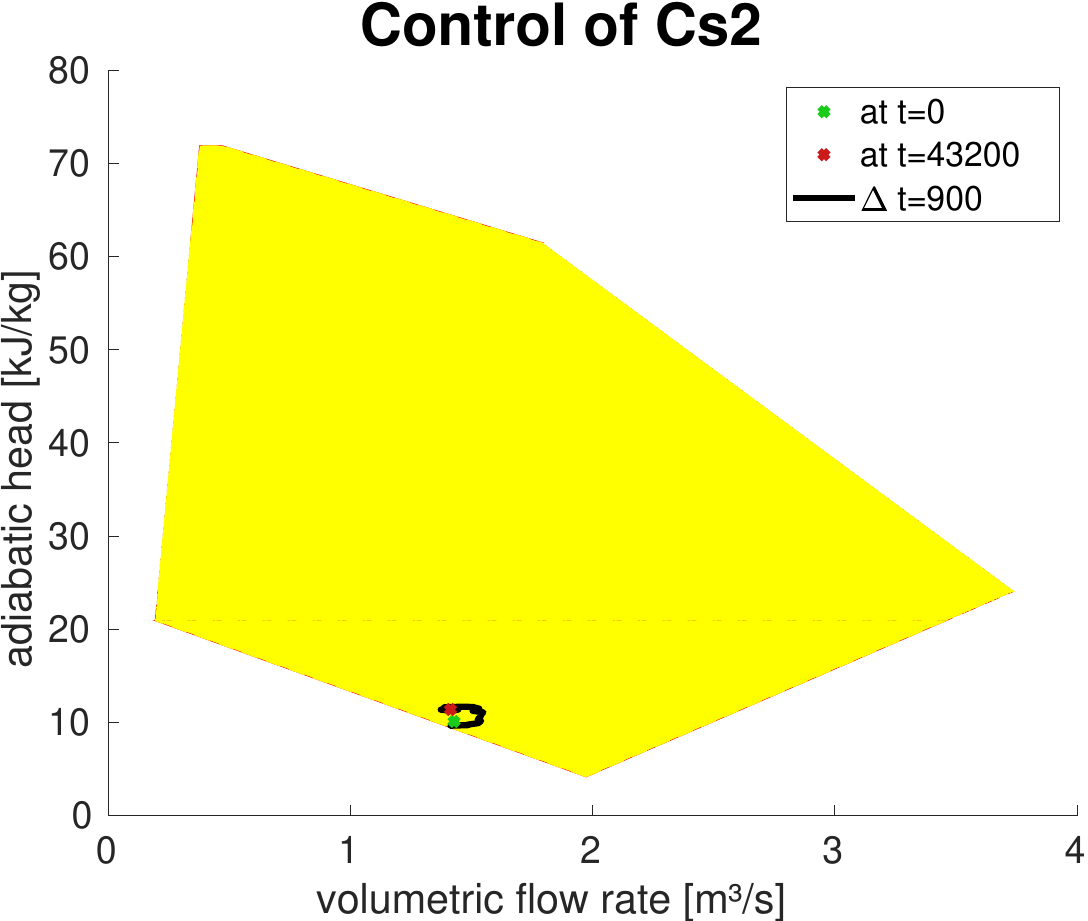}
    \end{subfigure}%
    ~
    \begin{subfigure}[t]{0.44\textwidth}
        \centering
        \includegraphics[width=\textwidth]{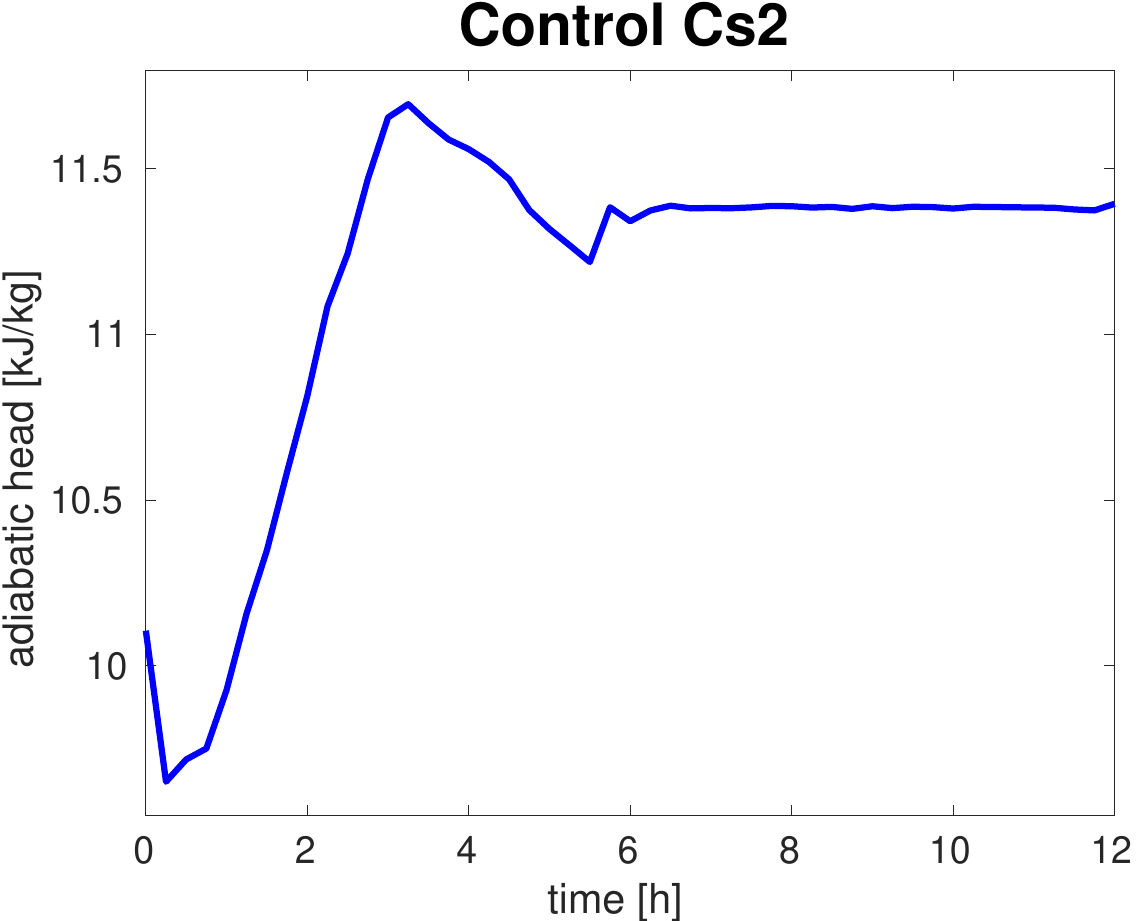}
    \end{subfigure}
    \\[0.5cm]
    \begin{subfigure}[t]{0.42\textwidth}
        \centering
        \includegraphics[width=\textwidth]{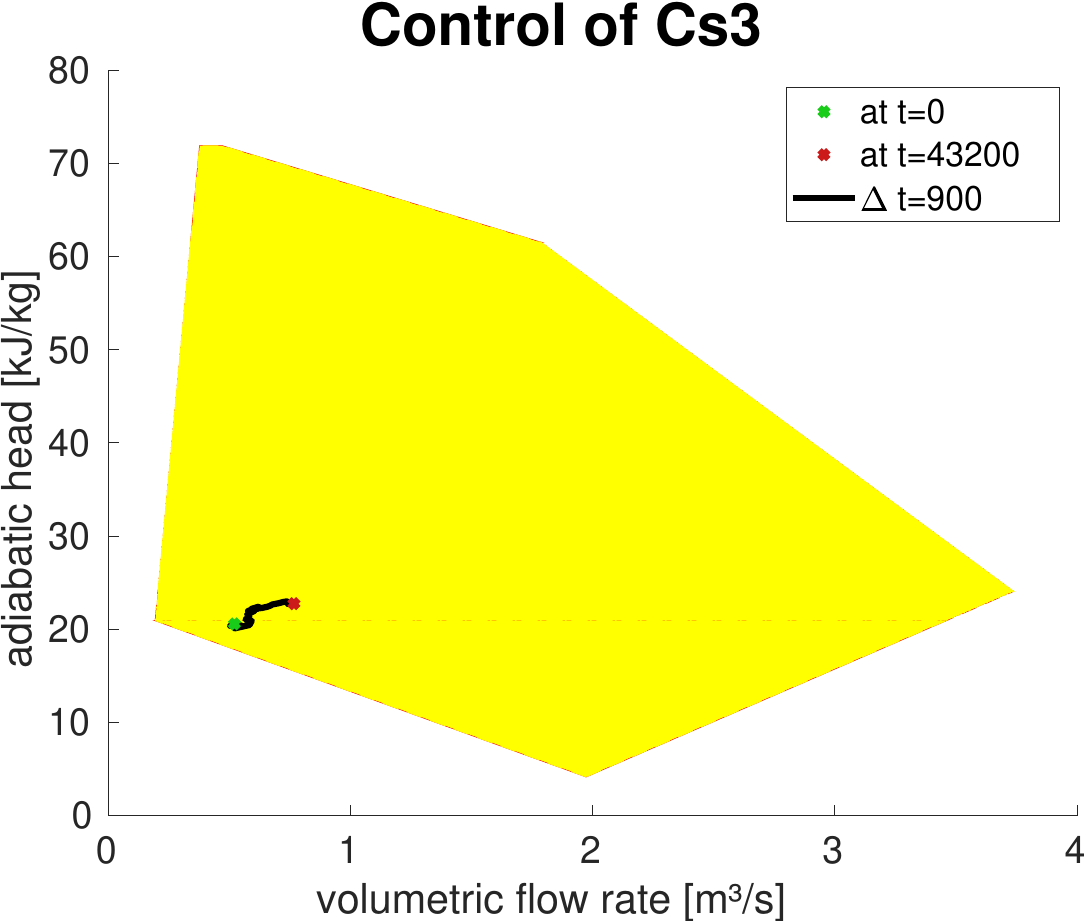}
    \end{subfigure}%
    ~
    \begin{subfigure}[t]{0.44\textwidth}
        \centering
        \includegraphics[width=\textwidth]{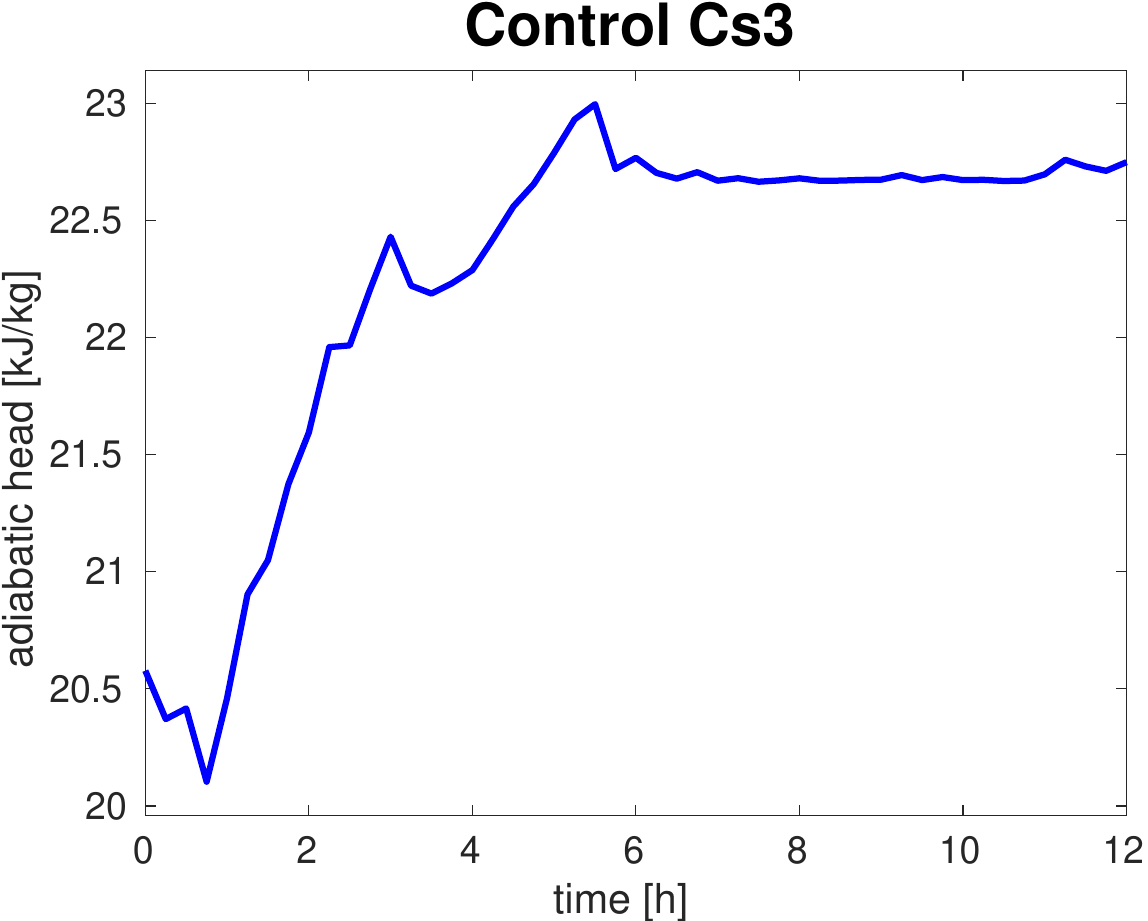}
    \end{subfigure}

    \caption{Controls $H_{ad}^{cp}(t)$ and constraints given by the
    $(Q,H_{ad})$-diagrams of the compressor stations cp=Cs1, Cs2, Cs3 (top to
    bottom)}
    \label{figure:controlsAndConstraintsSteinach}
\end{figure*}

At time $t=0$, the network is in a stationary state $A$, resulting from a simulation with stationary boundary conditions. The aim is to operate the network as to reach a stationary state $B$ under certain constraints. For the sinks T2 and T3, we have the following targets:
\begin{align*}
 q_{T2,target}(t) &= \begin{cases}
               \SI{40}{\meter\cubed\per\second},		& 0\leq t < \SI{7200}{\second}\,,\\
               \SI{5}{\meter\cubed\per\second\squared} \cdot t + \SI{30}{\meter\cubed\per\second},  & \SI{7200}{\second} \leq t < \SI{21600}{\second}\,,\\
               \SI{60}{\meter\cubed\per\second},		& \SI{21600}{\second} \leq t \leq \SI{43200}{\second}\,,
              \end{cases}\\
 q_{T3,target}(t) &= \begin{cases}
               \SI{60}{\meter\cubed\per\second},		& 0\leq t < \SI{7200}{\second}\,,\\
               \SI{-5}{\meter\cubed\per\second\squared} \cdot t + \SI{70}{\meter\cubed\per\second},  	& \SI{7200}{\second} \leq t < \SI{21600}{\second}\,,\\
               \SI{40}{\meter\cubed\per\second},		& \SI{21600}{\second} \leq t \leq \SI{43200}{\second}\,.
              \end{cases}
\end{align*}
They enter the objective function in \eqref{obj:Mnv} with weights $\beta_{T2}=10^{-4}$ and
$\beta_{T3}=10^{-5}$.
The target state B has to be accessed under certain pressure constraints at T2 and T3,
\begin{align*}
\text{T2:} & \quad p \geq p_{T2}^l(t) \quad\text{with}\;
p_{T2}^l(t) = \begin{cases}
				\SI{79}{\bar} + \SI{4/3}{\bar\per\hour} \cdot t & 0\leq t < \SI{3}{\hour}\,,\\
				\SI{83}{\bar} & \SI{3}{\hour} \leq t \leq \SI{12}{\hour}\,,
			\end{cases}\\
\text{T3:} & \quad p_{T3}^u\geq p \geq p_{T3}^l \quad\text{with}\;p_{T3}^u=\SI{61}{\bar},\;
p_{T3}^l=\SI{59}{\bar}.
\end{align*}
The constraints for the compressor stations Cs1-Cs3 are given by identical characteristic fields
as shown in Fig.~\ref{figure:controlsAndConstraintsSteinach}. We choose $\alpha=0.1$ as weight in
\eqref{obj:Mnv}.

\begin{figure}[!h]
    \centering
    \begin{subfigure}[t]{0.45\textwidth}
        \centering
        \includegraphics[width=\textwidth]{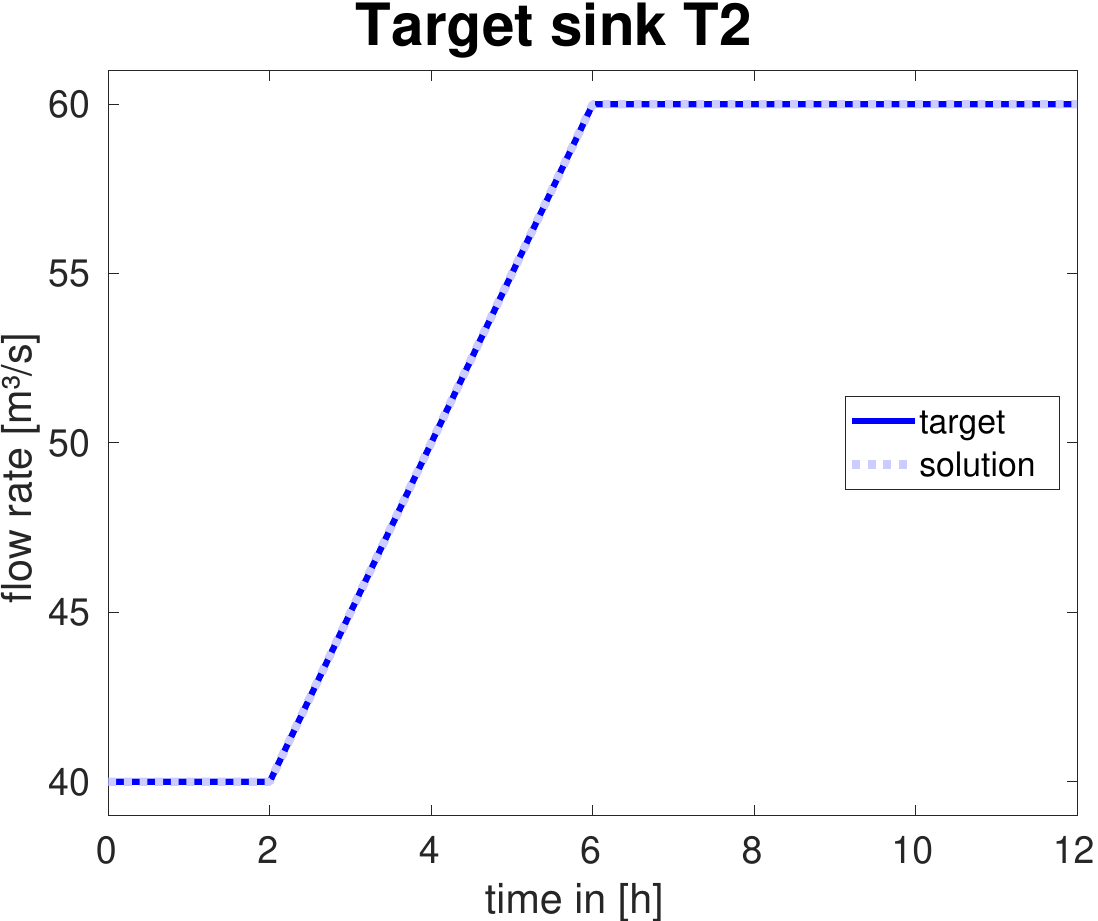}
    \end{subfigure}%
    ~
    \begin{subfigure}[t]{0.45\textwidth}
        \centering
        \includegraphics[width=\textwidth]{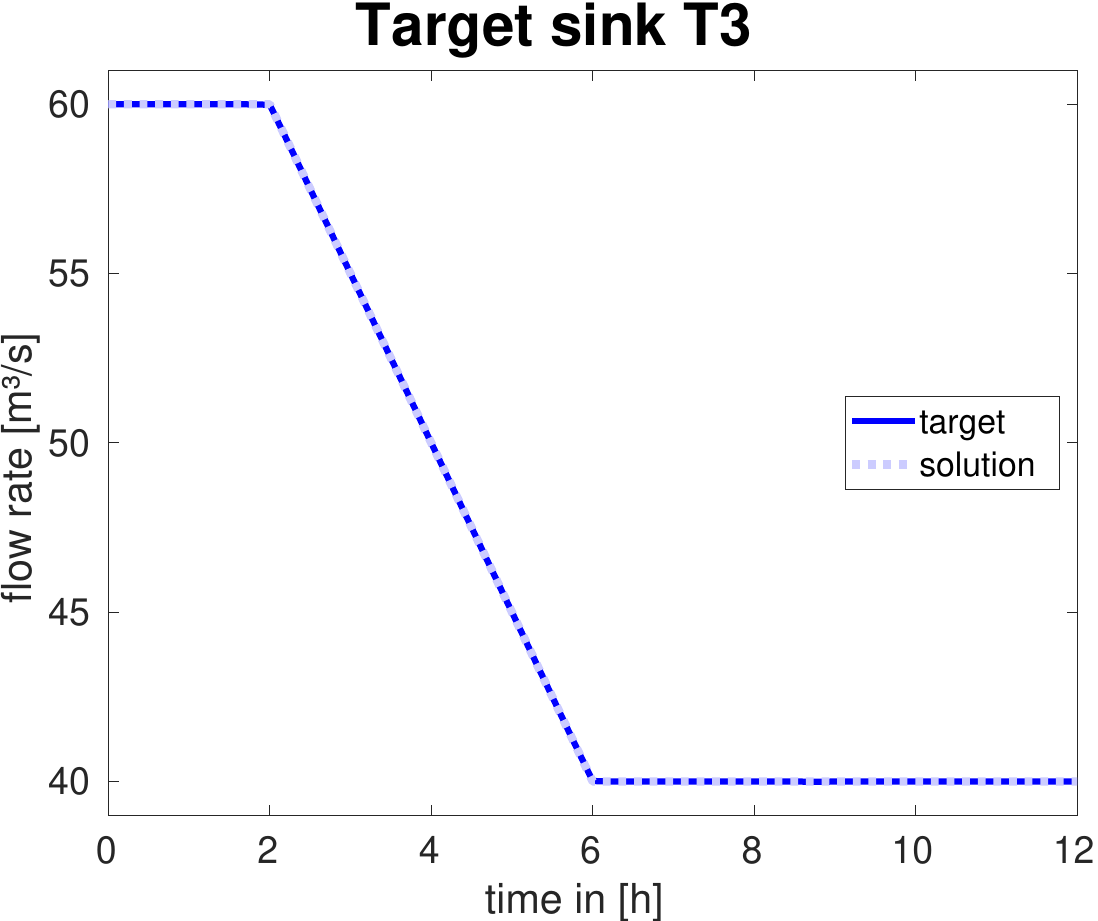}
    \end{subfigure}
    \\[0.5cm]
    \begin{subfigure}[t]{0.45\textwidth}
        \centering
        \includegraphics[width=\textwidth]{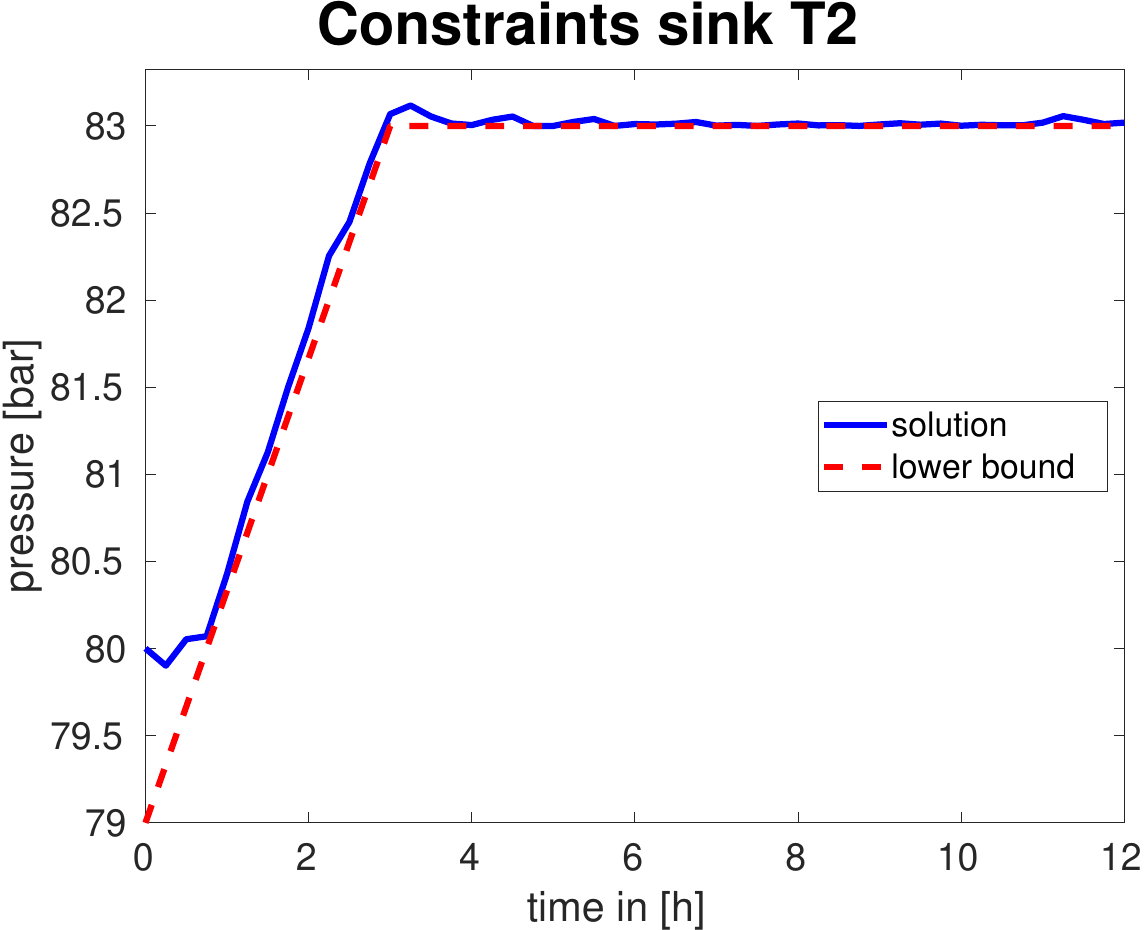}
    \end{subfigure}%
    ~
    \begin{subfigure}[t]{0.45\textwidth}
        \centering
        \includegraphics[width=\textwidth]{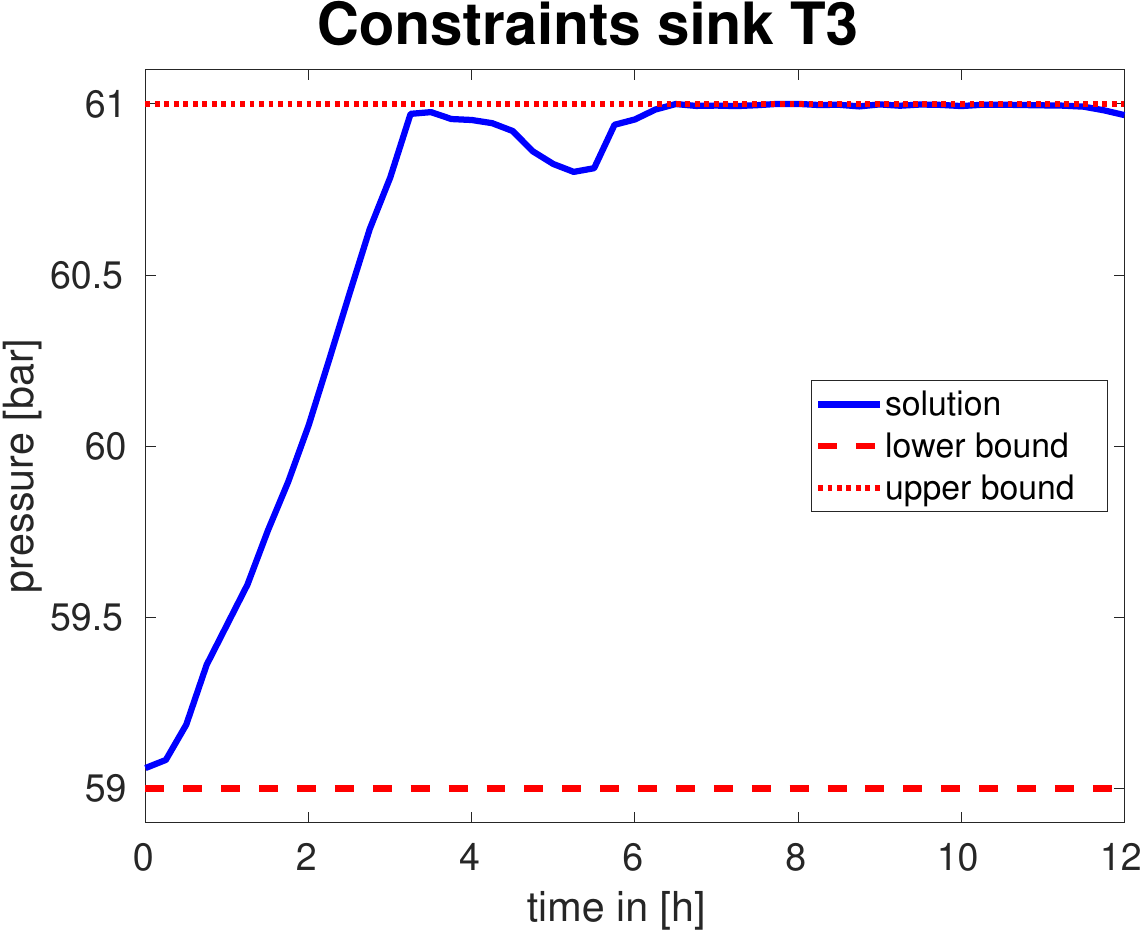}
    \end{subfigure}
    \caption{Controlled flow rate $q_s(t)$ (top) and pressure $p_s(t)$ (bottom)
    compared to target values $q_{s,target}(t)$ and constraints $p_s^{l,u}(t)$ for
    sinks s=T2,\,T3}
    \label{figure:targetAndConstraintSteinach}
\end{figure}

The simulation interval is [\SI{0}{\hour},\SI{12}{\hour}]. The discrete control vector $c\in\R^{240}$ consists
of the three adiabatic heads $H_{ad}^{Csk}(t_i),\,k=1,2,3,$ of the compressor stations and the flow rates $q_{Tk}(t_i),\,k=2,3,$
taken at the discrete control points $t_i=i\,\SI{900}{\second}$, $i=1,\ldots,48$. We use linear interpolation
between the control points whenever necessary. The optimization algorithm {\sc Donlp2} equipped with gradient information
delivers solutions for tolerance $\texttt{TOL}=\num{5e-4}$ after $16.5$ min computing time. The corresponding targets and constraints of the sinks T2 and T3 are shown in Figure~\ref{figure:targetAndConstraintSteinach}. We nicely see that the targets $q_{T2,target}(t)$ and $q_{T3,target}(t)$ are met very well by the control variables $q_{T2}(t)$ and $q_{T3}(t)$ during the optimization process and that also the constraints on
the pressures $p_{T2}(t)$ and $p_{T3}(t)$ at the sinks are maintained.

The compressor stations are controlled via their individual adiabatic heads.
These controls are shown in Figure~\ref{figure:controlsAndConstraintsSteinach}. On the left
hand side, the controls are drawn inside the characteristic field.
We see that the controls stay inside the characteristic field and thus meet the technical
constraints of the compressor stations.

Eventually, we would like to mention that relaxing the tolerance of {\sc Donlp2} to
$\texttt{TOL}=\num{5e-3}$ still yields satisfactory results after $7.5$ min computing time,
but the solutions found are less smooth.

\FloatBarrier
\subsection{The Greek Network (GasLib-134)}
\begin{figure}[!ht]
 \centering
 \includegraphics[width=0.35\textwidth]{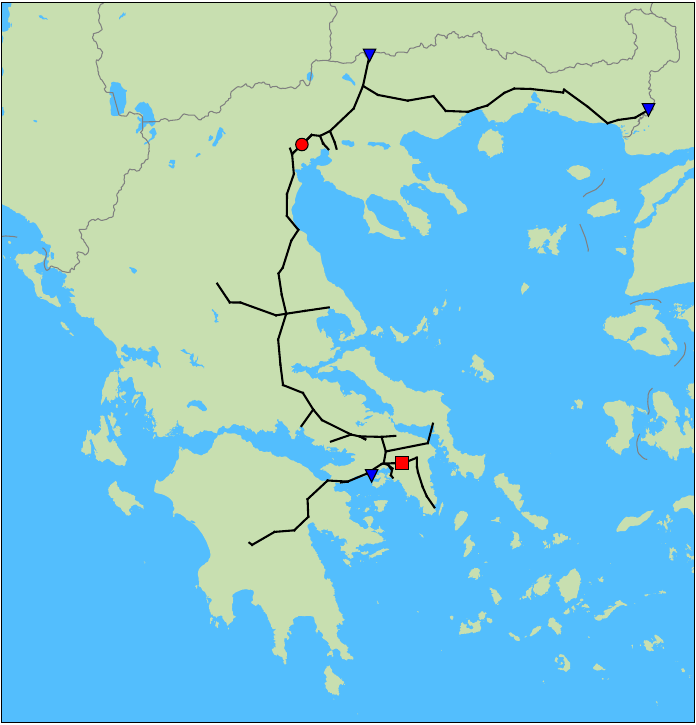}\hspace{1cm}
 \includegraphics[width=0.45\textwidth]{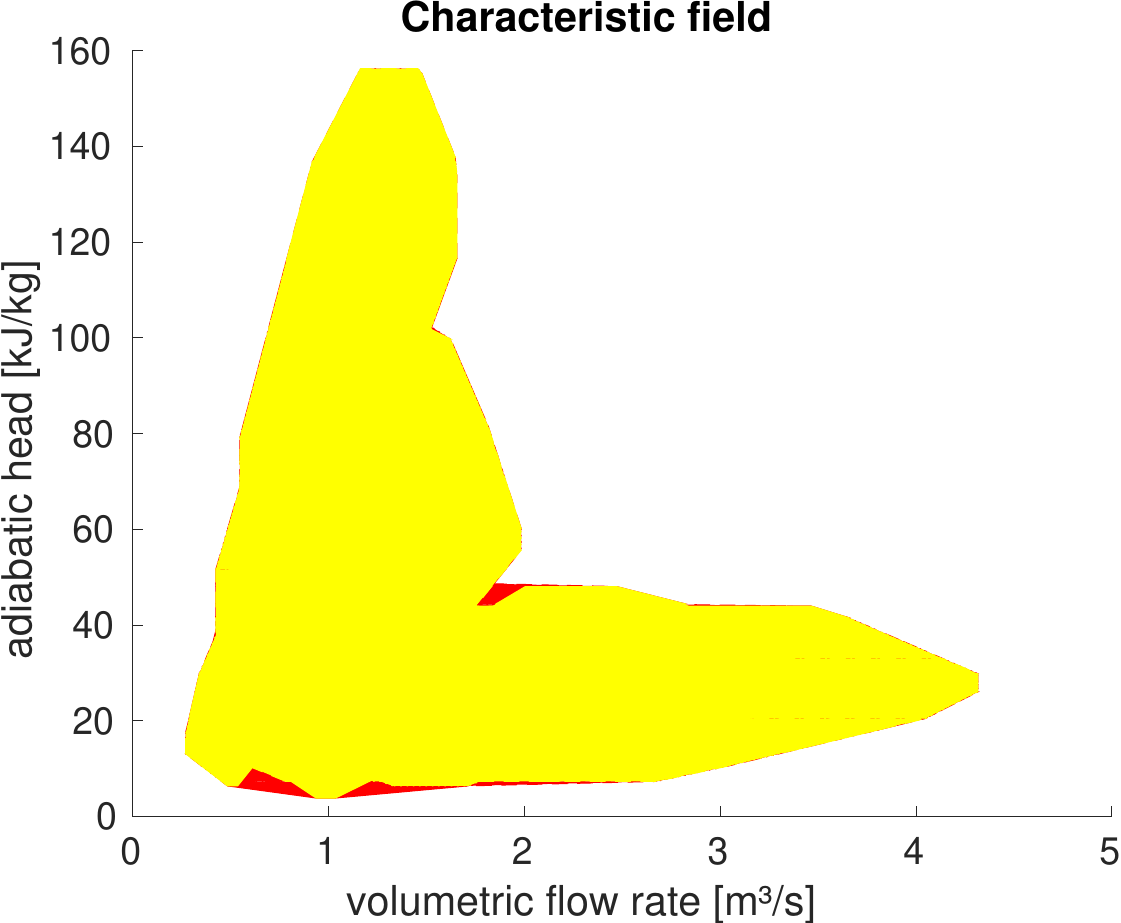}
 \caption{The Greek network (GasLib-134), reproduced from \url{gaslib.zib.de} (left),
 and characteristic field of the compressor station \texttt{cs} (right).}
 \label{figure:Griechenland}
\end{figure}

As second example, we consider the real Greek network with 134 nodes including 3~sources,
45~sinks as well as 133~edges, one compressor station and one control valve, see Fig.~\ref{figure:Griechenland}.
The specific data are available from \url{gaslib.zib.de}.

\begin{table}[h!]
\centering
\caption{Simulation results for the Greek network GasLib-134}
\begin{tabular}{ccccccc}
\hline
\TOL & rel. error & $M(u_h)$ & max/min $\Delta t$ & max/min $\Delta x$
& M1:M2:M3 [\%] & CPU [s] \\
\hline
1e-01 & 6.66e-03 & 0.0903765 & 7200/3600 & 18297.5/162.835 & 69:31:00 & 0.754 \\
1e-02 & 1.27e-03 & 0.0898926 & 7200/3600 & 18297.5/162.835 & 43:57:00 & 0.907 \\
1e-03 & 4.41e-04 & 0.0898186 & 1800/900 & 18297.5/162.835 & 14:84:02 & 3.387 \\
1e-04 & 6.71e-05 & 0.0897850 & 225/225 & 18297.5/162.835 & 08:81:11 & 14.74 \\
1e-05 & 2.36e-06 & 0.0897788 & 225/14.0625 & 18297.5/162.835 & 04:73:23 & 104.7 \\
\hline
\multicolumn{2}{c}{reference solution}  & 0.0897790 & 10 & 250/162.835 & 0:0:100 & 1273 \\
\hline
\end{tabular}
\label{table:GL134_SimulationResults}
\end{table}

First, we check typical simulation results for a given control for the compressor
station and varying tolerances. Exemplary results for $T=\SI{24}{\hour}$ are presented in Tab.~\ref{table:GL134_SimulationResults}.
Practically sufficient accuracy is achieved with $\TOL=10^{-2}$.

Second,
we are again interested in computing an optimized way to transfer the
gas network from a stationary state $A$ to another stationary state $B$ under
certain constraints. Except one, all sinks and sources have the same structure
for the target functional:
\begin{align*}
 q_{s,target}(t) &= \begin{cases}
               q_s^A,	& 0\leq t < t_1\,,\\
               q_s^A + \frac{t-t_1}{t_2-t_1}\left(q_s^B-q_s^A\right),  & t_1 \leq t < t_2\,,\\
               q_s^B,	& t_2 \leq t \leq T\,.
           \end{cases}
\end{align*}
for $s=1,\ldots,47$, and with $t_1=\SI{2}{\hour}$, $t_2=\SI{6}{\hour}$, and simulation time $T=\SI{24}{\hour}$.
At source \verb|node_20|, which can be identified to be the small blue triangle in the very north of Greece in Figure~\ref{figure:Griechenland}, we prescribe the pressure in order to get a well-defined stationary state. We set
\begin{align*}
 p_{20}(t) &= \begin{cases}
               p_{20}^A	& 0\leq t < t_1\,,\\
               p_{20}^A + \frac{t-t_1}{t_2-t_1}\left(p_{20}^B-p_{20}^A\right)  & t_1 \leq t < t_2\,,\\
               p_{20}^B	& t_2 \leq t \leq T\,.
           \end{cases}
\end{align*}
\begin{figure*}[t!]
    \centering
    \begin{subfigure}[t]{0.42\textwidth}
        \centering
        \includegraphics[width=\textwidth]{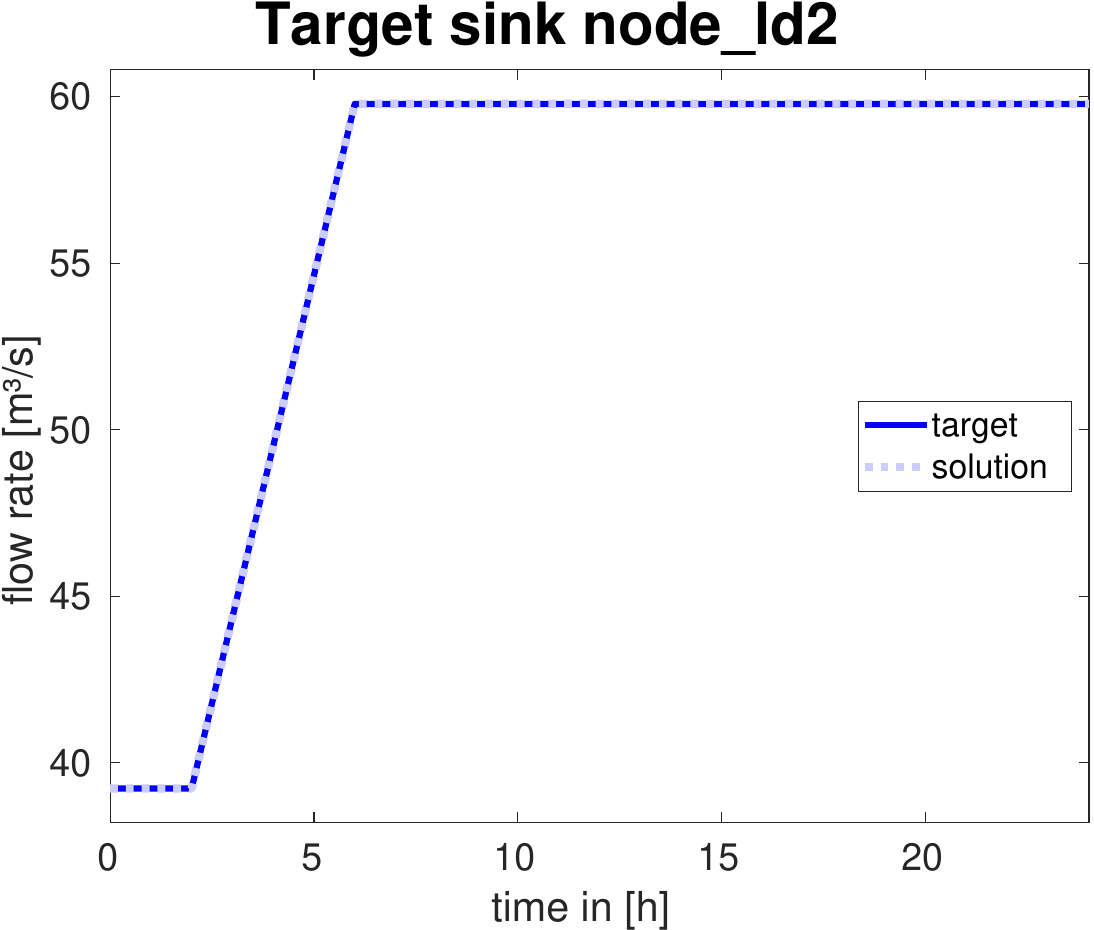}
    \end{subfigure}%
    ~
    \begin{subfigure}[t]{0.42\textwidth}
        \centering
        \includegraphics[width=\textwidth]{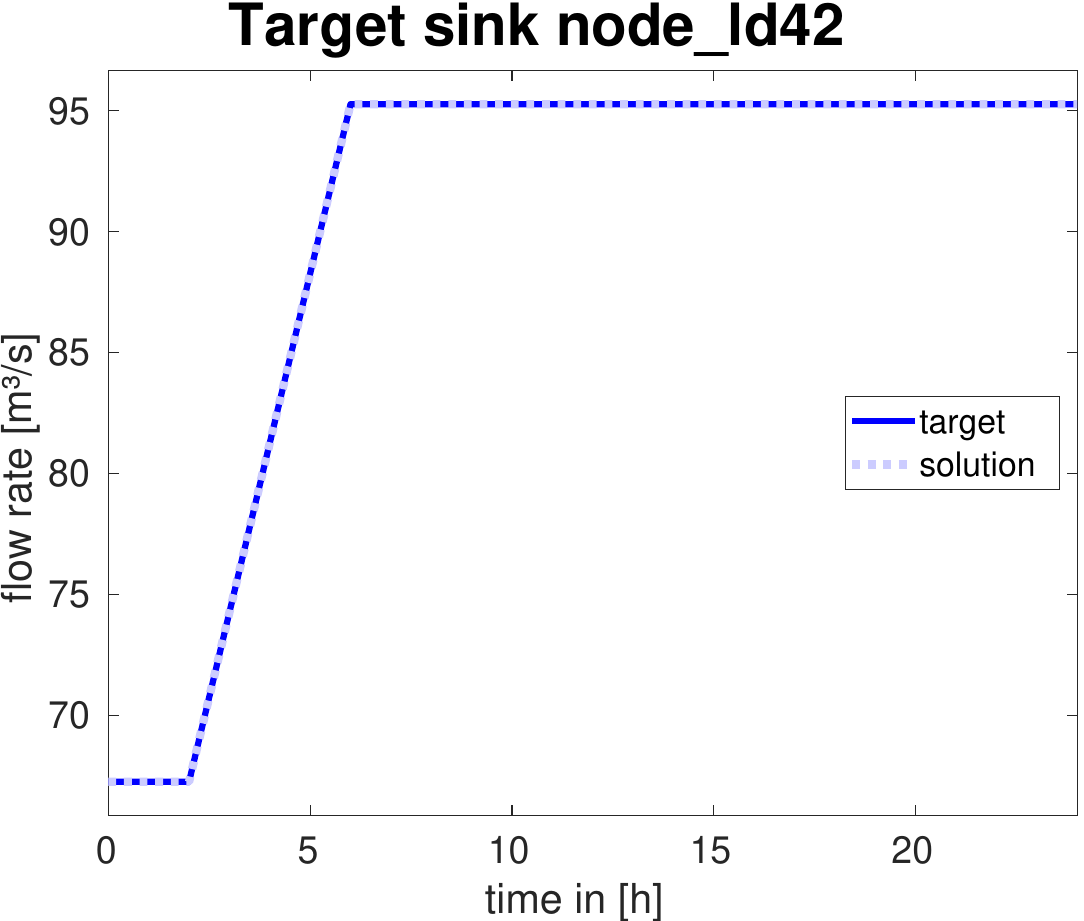}
    \end{subfigure}
    \\[0.5cm]
    \begin{subfigure}[t]{0.42\textwidth}
        \centering
        \includegraphics[width=\textwidth]{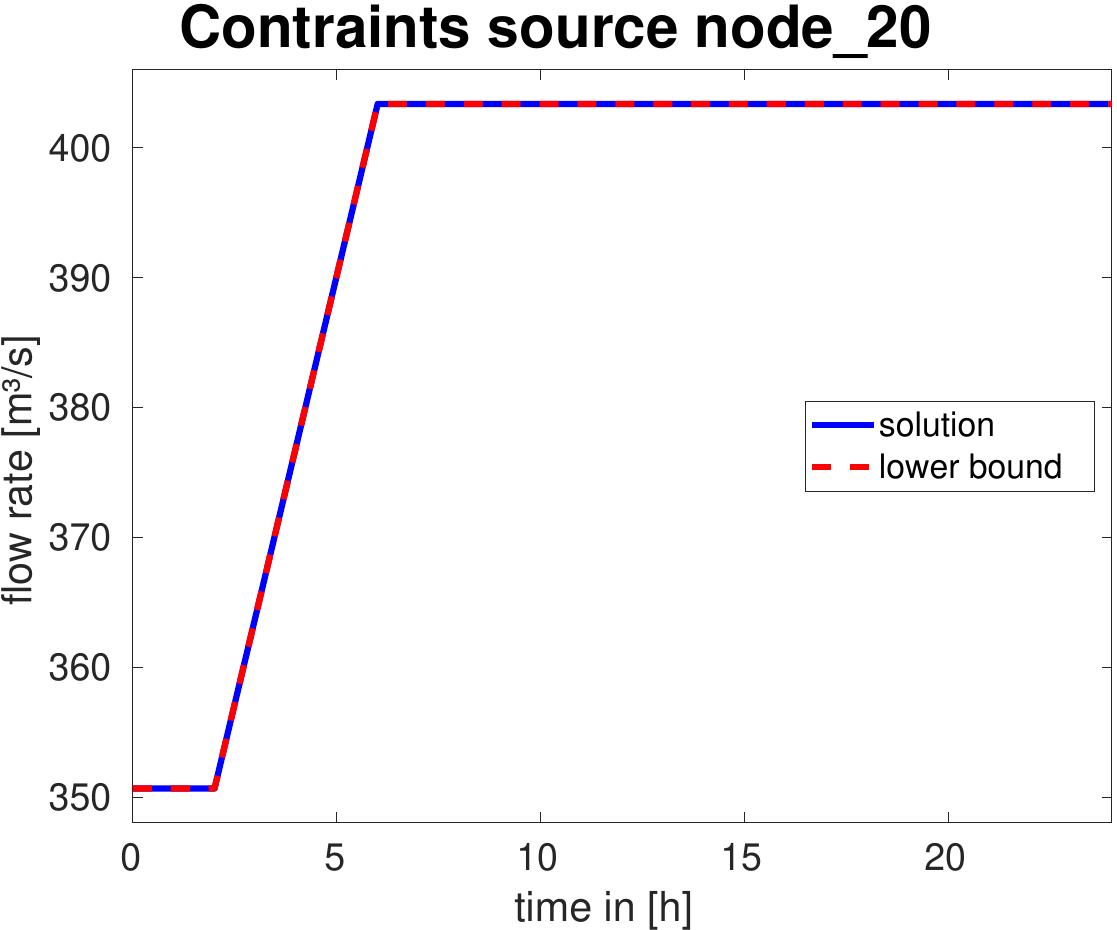}
    \end{subfigure}%
    ~
    \begin{subfigure}[t]{0.42\textwidth}
        \centering
        \includegraphics[width=\textwidth]{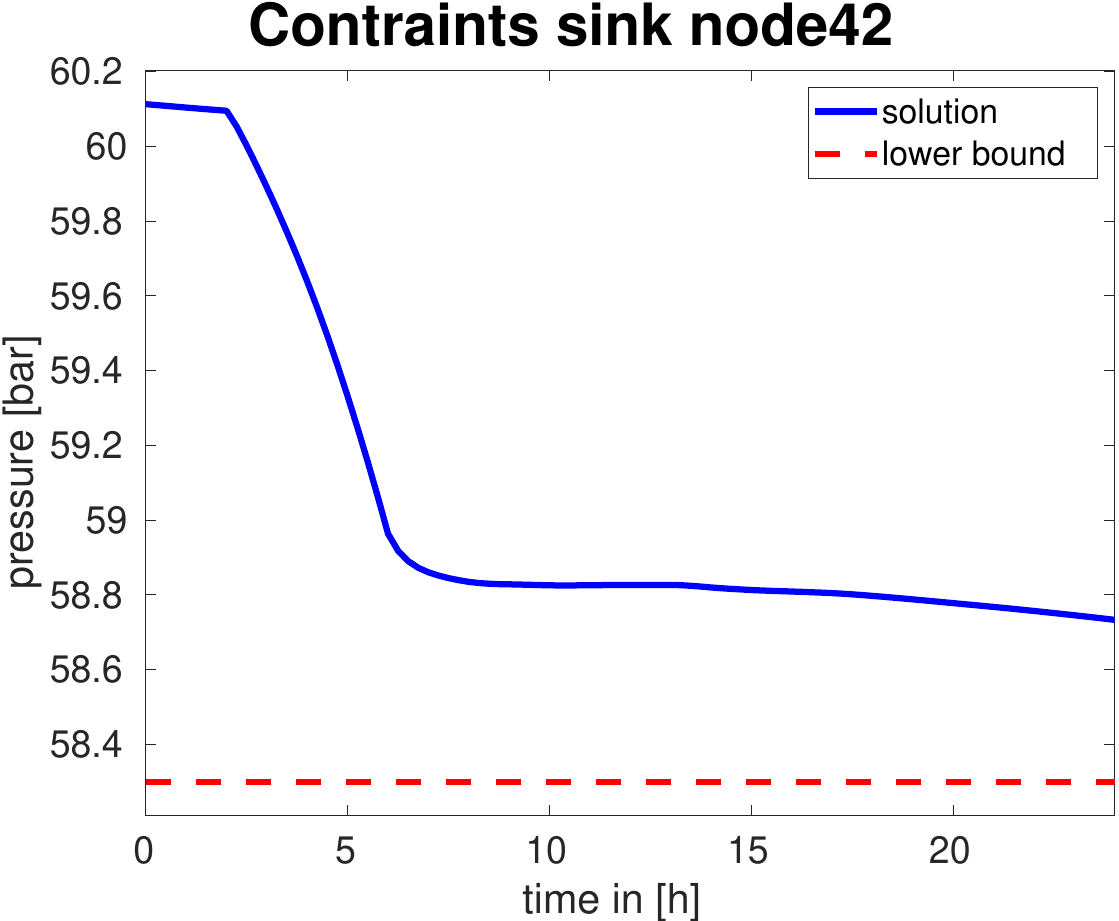}
    \end{subfigure}
    \caption{Controlled flow rates $q_s(t)$ compared to target values $q_{s,target}(t)$
    for $s=2,42$ (top), flow rate $q_{20}(t)$ and pressure $p_{42}(t)$
    compared to constraints $q_{20}^l$ and $p_{42}^l$ (bottom) for discrete control points
    with $\Delta t=\SI{900}{\second}$}
    \label{figure:targetAndConstraintGreek}
\end{figure*}
In order to achieve the stationary state $B$ in the end, we set the lower bound
for the flow rate at source \verb|node_20| to be
\begin{align*}
 q_{20}^l(t) &=\begin{cases}
               q_{20}^A	& 0\leq t < t_1\,,\\
               q_{20}^A + \frac{t-t_1}{t_2-t_1}\left(q_{20}^B-q_{20}^A\right)  & t_1 \leq t < t_2\,,\\
               q_{20}^B	& t_2 \leq t \leq T\,.
           \end{cases}
\end{align*}
Further constraints are defined for the sinks \verb|node_42| and \verb|node_45|,
\begin{align*}
p_{42}(t)\geq p_{42}^l=\SI{58.3}{\bar} \quad\text{and}\quad p_{45}(t)\geq p_{45}^l=\SI{61}{\bar}.
\end{align*}
They are located in the very south and south-west of the network, respectively.
The constraints for the compressor station \verb|cs| is given by a characteristic field as shown in Fig.~\ref{figure:Griechenland}. We employ the semiconvex approximation of the polyhedric approximation as described above, see Fig.~\ref{abb:aggKennfeld}.

\begin{figure*}
   \centering
    \begin{subfigure}[t]{0.42\textwidth}
        \centering
        \includegraphics[width=\textwidth]{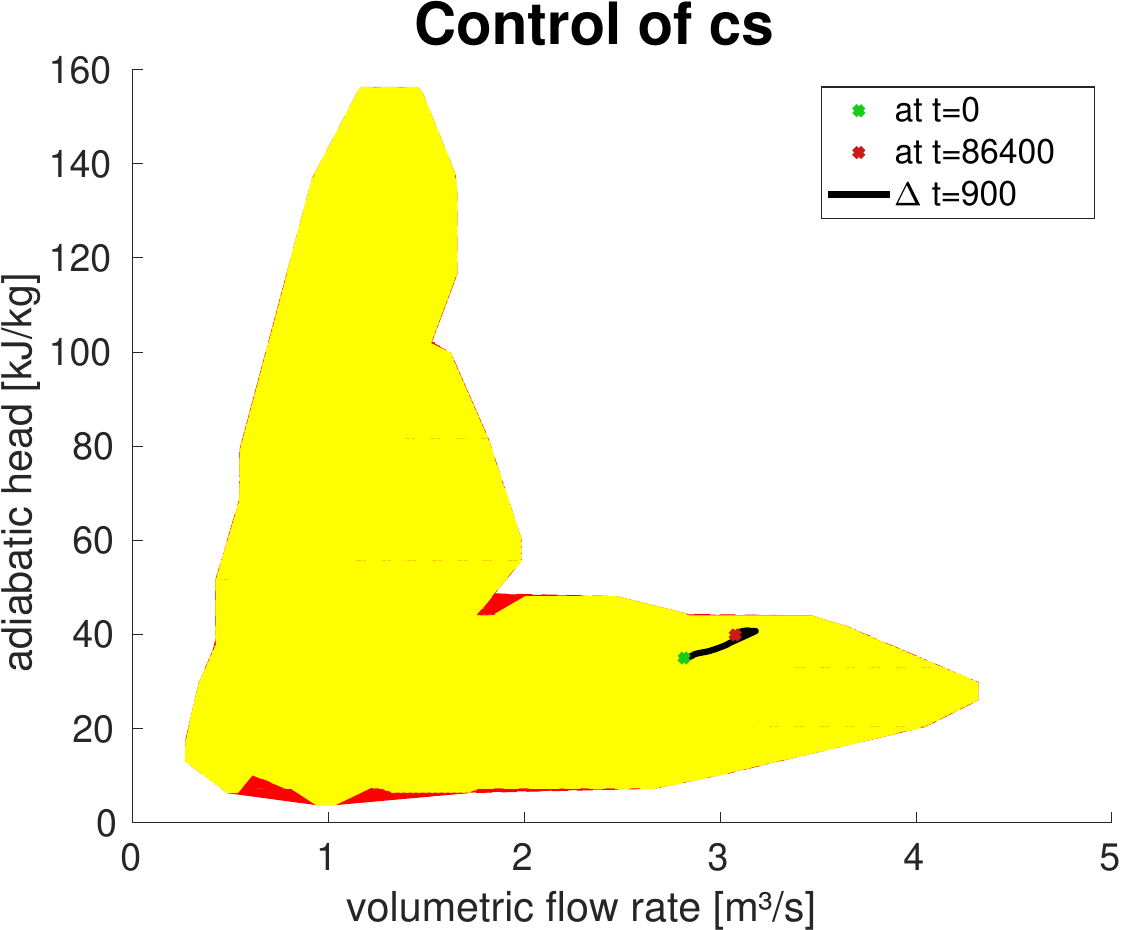}
    \end{subfigure}%
    ~
    \begin{subfigure}[t]{0.42\textwidth}
        \centering
        \includegraphics[width=\textwidth]{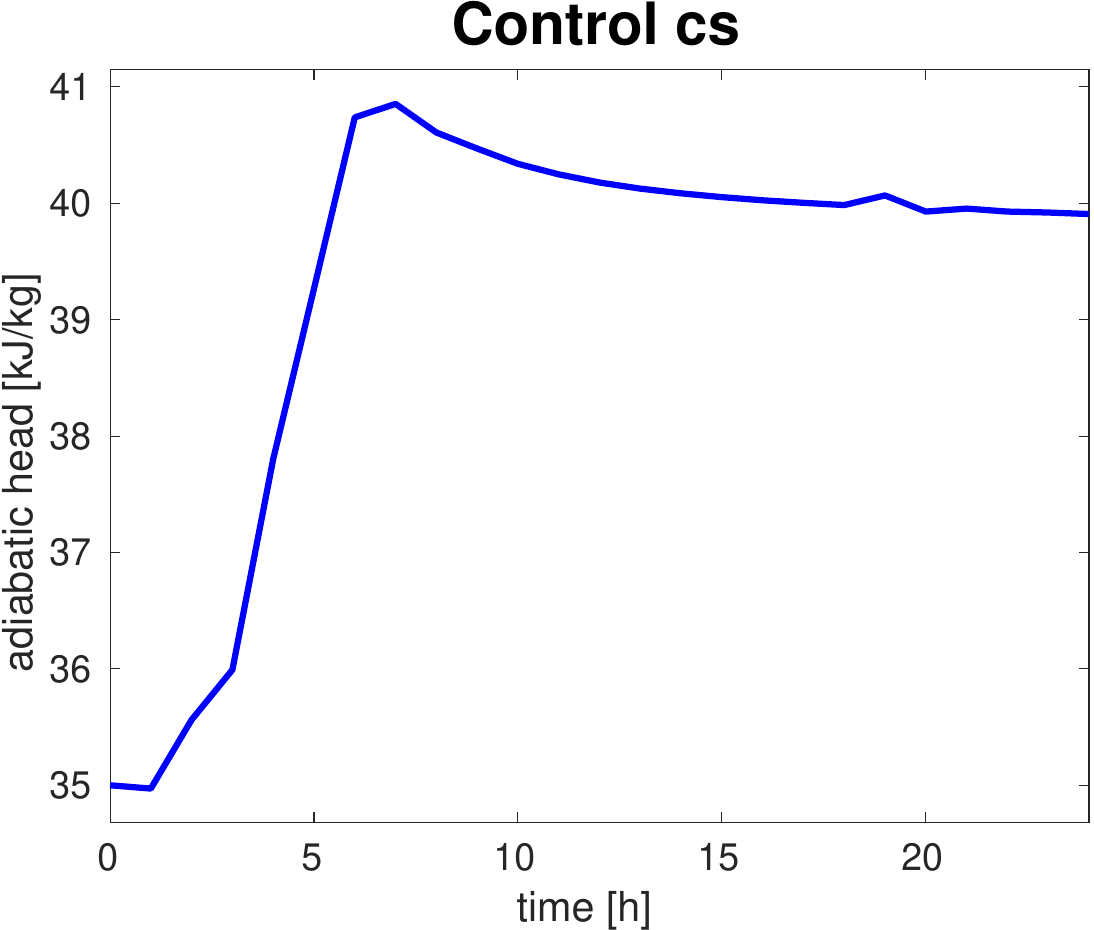}
    \end{subfigure}
      \\[0.5cm]
    \begin{subfigure}[t]{0.42\textwidth}
        \centering
        \includegraphics[width=\textwidth]{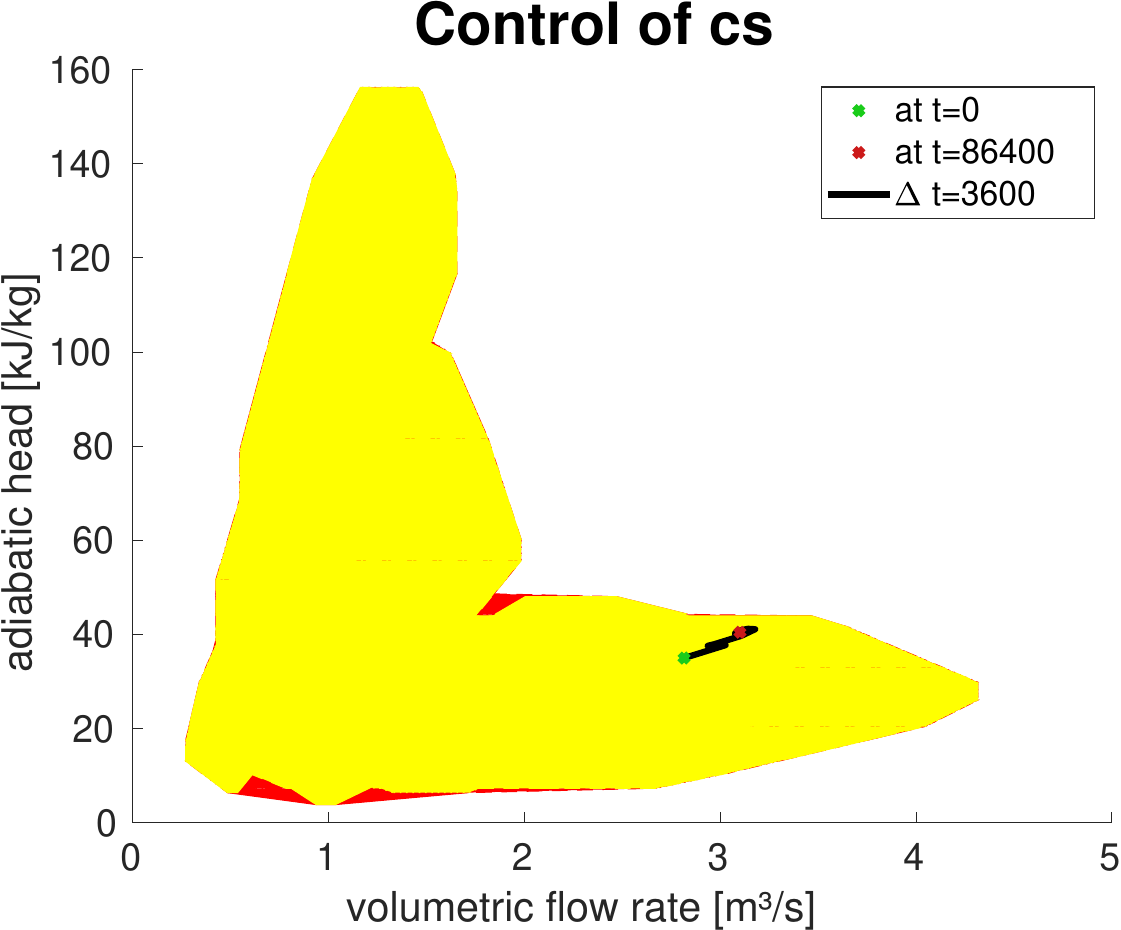}
    \end{subfigure}%
    ~
    \begin{subfigure}[t]{0.42\textwidth}
        \centering
        \includegraphics[width=\textwidth]{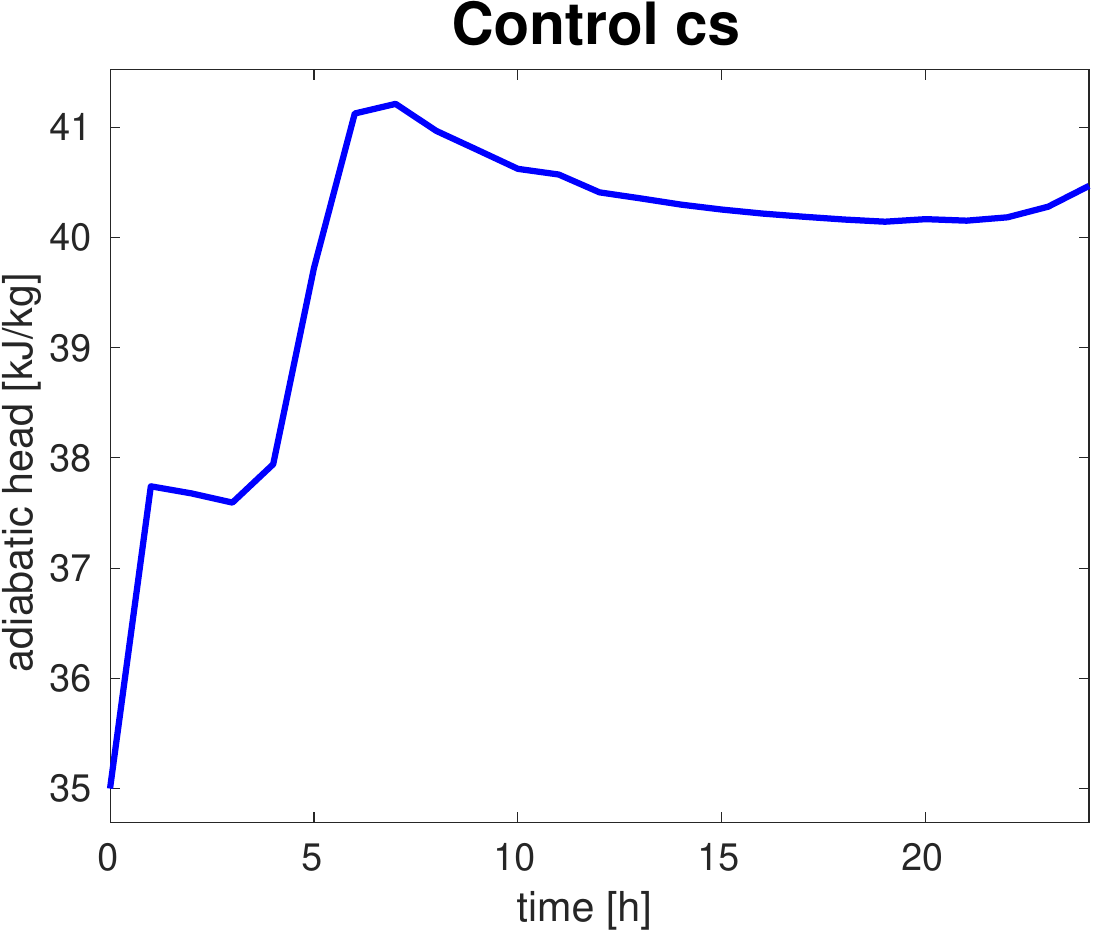}
    \end{subfigure}
    \caption{Control $H_{ad}^{cs}(t)$ and constraints given by the
    $(Q,H_{ad})$-diagram of the compressor station \texttt{cs} for discrete control points
    with $\Delta t=\SI{900}{\second}$ (top) and $\Delta t=\SI{3600}{\second}$ (bottom)}
    \label{figure:controlsAndConstraintsGreek}
\end{figure*}

The simulation interval is now [\SI{0}{\hour},\SI{24}{\hour}].
Beside the adiabatic head $H_{ad}^{cs}(t_i)$ for the
compressor station, we consider all flow rates $q_s(t_i)$, $s=1,\ldots,47$, at the discrete control points
$t_i=i\,\SI{900}{\second}$, $i=1,\ldots,96$ as control variables, which gives $c\in\R^{4608}$. As a
consequence, we enforce $\Delta t\le \SI{900}{\second}$ in our simulations to reach every control
point with our adaptive scheme. We further set
$\beta_s=10^{-5}$ for all $s$ and $\alpha=1$ in the objective function \eqref{obj:Mnv}.
The optimization is challenging due to the additional state constraints.
It was run with {\sc Donlp2} with a tolerance $\texttt{TOL}=\num{e-3}$ and took 10.1~min. We note that
the dimension of the overall discrete state vector taken over space and time reaches more than a million in
certain cases.
Fig.~\ref{figure:targetAndConstraintGreek} shows exemplary results for the targets of two selected nodes.
The variable operation mode for the compressor station is plotted in
Fig.~\ref{figure:controlsAndConstraintsGreek} (top). All constraints are satisfied.

We have also reduced the discrete control points to $t_i=i\,\SI{3600}{\second}$, $i=1,\ldots,24$, which scales
down the dimension of the control space to $c\in\R^{1152}$. {\sc Donlp2} was still able to find a
feasible solution after 2.5~min. The corresponding operation mode for the compressor station is shown in
Fig.~\ref{figure:controlsAndConstraintsGreek} (bottom). It differs from the solution found for
$\Delta t=\SI{900}{\second}$.

\section{Conclusion and Outlook} \label{sec:conclusion}
We have presented a novel optimization approach to support a transient management
of large-scale gas networks with real semiconvex models for aggregated
characteristic fields of compressor stations.
Given the significantly reduced computational time for a one-time simulation by
fully self-adaptive space-time-model discretizations and the associated inherent error control,
the method proposed is a promising candidate for a reliable work horse in a predictive transient
software framework that creates gas flow schedules and forecasts in near real-time. As
an example, we
have studied the practically important issue of nomination validations, i.e., transferring
the gas network from a stationary state to another stationary state usually described by a
change in the gas demands by consumers. During this process certain state constraints have
to be satisfied, making the corresponding optimal control problem challenging. We could
observe that treating the new outflow rates as control variables and adding regularization terms
of tracking type for them to the objective function is an attractive way to find answers
to the question whether a desired steady state can be reached.

In future projects, we will include mixed integer formulations to represent network
topology changes and use probabilistic constraint optimization to investigate the
influence of uncertainties.


\section*{Declarations}
\textit{Funding}. The authors P.~Domschke and J.~Lang are supported by the Deutsche Forschungsgemeinschaft (DFG, German Research Foundation) within the collaborative research center TRR154 “Mathematical modeling, simulation and optimisation using the example of gas networks“ (Project-ID239904186, TRR154/2-2018, TP B01).\\

\noindent\textit{Conflicts of interest/Competing interests}. The authors declare that they have no known competing financial interests or personal relationships that could have appeared to influence the work reported in this paper.\\

\noindent\textit{Availability of data and material}. Detailed descriptions of the gas networks GasLib-582 and
GasLib-134 can be found at gaslib.zib.de.\\

\noindent\textit{Code availability}. The nonlinear programming package {\sc Donlp2} by Peter Spellucci is free for academic use only, commercial use needs licensing.

\bibliographystyle{plain}
\bibliography{gasliterature}

\end{document}